\documentclass[12pt]{article}

\usepackage[margin=1in]{geometry} 
\usepackage{graphicx}            
\usepackage{amsfonts}  
\usepackage[breaklinks]{hyperref}
\usepackage{tikz}
\usepackage{float}
\usetikzlibrary{patterns,arrows,decorations.pathreplacing}
\usetikzlibrary{decorations.markings}
\usetikzlibrary{calc,shapes}

\newcommand{\s}[1]{\left\lvert #1 \right\rvert}
\newcommand{\br}[1]{\left( #1 \right)}
\newcommand{\J}{\widetilde{\mathcal{J}}}
\newcommand{\eps}{\varepsilon}
\newcommand{\m}{\widetilde{m}}
\newcommand{\R}{\widetilde{\mathcal{R}}}
\newcommand{\tF}{\widetilde{F}}
\newcommand{\noth}{\varnothing}

\newcommand{\tY}{\widetilde{\mathcal{Y}}}

\usepackage{amsmath,amsthm,amssymb,enumerate,mathrsfs,mathtools}
\usepackage{relsize}
\usepackage{a4wide}
\usepackage{verbatim}
\usepackage{xcolor}

\usepackage{makecell}
\usepackage[numbers,sort&compress]{natbib}

\usepackage{enumitem}
\usepackage[capitalise]{cleveref}
\crefname{equation}{}{}
\crefname{prop}{Proposition}{Propositions}

\usepackage[utf8]{inputenc}

\newtheorem{thm}{Theorem}[section]
\crefname{thm}{Theorem}{Theorems}
\newtheorem{lem}[thm]{Lemma}
\newtheorem{prop}[thm]{Proposition}
\newtheorem{cor}[thm]{Corollary}

\newtheorem{ex}[thm]{Example}

\newtheorem{claim}[thm]{Claim}
\newtheorem{fact}[thm]{Fact}

\theoremstyle{definition}
\newtheorem{defn}[thm]{Definition}

\numberwithin{equation}{section}

\DeclareMathOperator{\blue}{blue}
\DeclareMathOperator{\red}{red}
\def\al#1{}
	\renewcommand{\al}[1]{\footnote{\textbf{AL: }#1}}  
    
\newenvironment{marknew}{\color{red} }{ }
\newcommand{\mn}{\begin{marknew}}
\newcommand{\umn}{\end{marknew}}

\newenvironment{proofclaim}[1][Proof of Claim]{\begin{proof}[#1]}{\end{proof}}

\usepackage[UKenglish]{babel}
\usepackage[UKenglish]{isodate}

\newcommand{\EMAIL}[1]{\textit{{E-mail}}: \texttt{\href{mailto:#1}{#1}}} 

\title{Towards Lehel's conjecture for $4$-uniform tight cycles}

\author{
Allan Lo 
\thanks{School of Mathematics, University of Birmingham, \EMAIL{s.a.lo@bham.ac.uk}} 
\and Vincent Pfenninger \thanks{School of Mathematics, University of Birmingham, \EMAIL{v.pfenninger@bham.ac.uk}}
}

\date{\today}

\begin{document}

\maketitle

\begin{abstract}
A $k$-uniform tight cycle is a $k$-uniform hypergraph with a cyclic ordering of its vertices such that its edges are all the sets of size~$k$ formed by~$k$ consecutive vertices in the ordering.
We prove that every red-blue edge-coloured~$K_n^{(4)}$ contains a red and a blue tight cycle that are vertex-disjoint and together cover~$n-o(n)$ vertices. Moreover, we prove that every red-blue edge-coloured~$K_n^{(5)}$ contains four monochromatic tight cycles that are vertex-disjoint and together cover~$n-o(n)$ vertices.
\end{abstract}

\section{Introduction}
An \emph{$r$-edge-colouring} of a graph (or hypergraph) is a colouring of its edges with~$r$ colours. A \emph{monochromatic subgraph} of an $r$-edge-coloured graph is one in which all the edges have the same colour.

Lehel conjectured that every $2$-edge-colouring of the complete graph on~$n$ vertices admits a partition of the vertex set into two monochromatic cycles of distinct colours, where the empty set, a single vertex and a single edge are considered to be \emph{degenerate cycles}. 
This conjecture was proved for large~$n$ by Łuczak, Rödl and Szemer{\'e}di~\cite{Luczak1998} using Szemer\'edi's Regularity Lemma. Allen~\cite{Allen2008} improved the bound on~$n$ by giving a different proof. 
Finally Bessy and Thomass\'e~\cite{Bessy2010} proved Lehel's conjecture for all $n \geq 1$.

Similar problems have also been considered for colourings with a general number of colours. In particular, a lot of attention has been given to the problem of determining the number of monochromatic cycles that are needed to partition an $r$-edge-coloured complete graph.
Erd\H{o}s, Gy\'arf\'as and Pyber~\cite{Erdos1991} proved that every $r$-edge-coloured complete graph can be partitioned into $O(r^2\log r)$ monochromatic cycles and conjectured that~$r$ monochromatic cycles would suffice. 
Their result was improved by Gy\'arf\'as, Ruszink\'o, S\'ark\"ozy and Szemer\'edi~\cite{Gyarfas2006} who showed that $O(r\log r)$ monochromatic cycles are enough. 
However, Pokrovskiy~\cite{Pokrovskiy2014} disproved the conjecture and proposed a weaker version of the conjecture that each $r$-edge-coloured complete graph contains~$r$ monochromatic vertex-disjoint cycles that together cover all but at most~$c_r$ of the vertices, where~$c_r$ is a constant depending only on~$r$. Pokrovskiy \cite{Pokrovskiy2016} subsequently proved that we can take $c_3 \leq 43000$ for large enough~$n$. 

Recently, generalisations of Lehel's conjecture to hypergraphs have also been considered.
For any positive integer~$k$, a \emph{$k$-uniform hypergraph}, or \emph{$k$-graph},~$H$ is an ordered pair of sets~$(V(H),E(H))$ such that $E(H) \subseteq \binom{V(H)}{k}$, where $\binom{S}{k}$ is the set of all subsets of~$S$ of size~$k$.
We abuse notation by identifying the $k$-graph~$H$ with its edge set~$E(H)$.
Hence by~$\s{H}$ we mean the number of edges of~$H$.
Let~$K_n^{(k)}$ be the complete $k$-graph on~$n$ vertices.

In $k$-graphs there are several notions of cycle. 
For integers $1 \leq \ell < k < n$, a $k$-graph~$C$ on~$n$ vertices is called an \emph{$\ell$-cycle} if there is an ordering of its vertices $V(C) = \{v_0, \dots, v_{n-1}\}$ such that $E(C) = \{ \{v_{i(k-\ell)},\dots, v_{i(k-\ell)+k-1}\}\colon 0 \leq i \leq n/(k-\ell) -1\}$, where the indices are taken modulo~$n$. 
That is, an $\ell$-cycle is a $k$-graph with a cyclic ordering of its vertices such that its edges are sets of~$k$ consecutive vertices and consecutive edges share exactly~$\ell$ vertices. 
(Note that~$k-\ell$ divides~$n$.)
A single edge or any set of fewer than~$k$ vertices is considered to be a \emph{degenerate $\ell$-cycle}.
Further, $1$-cycles and $(k-1)$-cycles are called \emph{loose cycles} and \emph{tight cycles}, respectively.

For loose cycles, Gy\'arf\'as and S\'ark\"ozy~\cite{gyarfas2013} showed that every $r$-edge-coloured complete $k$-graph on~$n$ vertices can be partitioned into~$c(k,r)$ monochromatic loose cycles. 
S\'ark\"ozy~\cite{Sarkozy2014} showed that, for~$n$ sufficiently large, $50kr\log(kr)$ loose cycles are enough.
For tight cycles, Bustamante, Corsten, Frankl, Pokrovskiy and Skokan~\cite{Bustamante2020} showed that every $r$-edge-coloured complete $k$-graph can be partitioned into~$C(k,r)$ monochromatic tight cycles. 
See \cite{Gyarfas2016} for a survey on other results about monochromatic cycle partitions and related problems.

In this paper, we investigate monochromatic tight cycle partitions in $2$-edge-coloured complete $k$-graphs on~$n$ vertices.
When $k=3$, Bustamante, H\`an and Stein \cite{Bustamante2017} showed that there exist two vertex-disjoint monochromatic tight cycles of distinct colours covering all but at most~$o(n)$ of the vertices.
Recently, Garbe, Mycroft, Lang, Lo and Sanhueza-Matamala \cite{Garbe2019} proved that two monochromatic tight cycles are sufficient to cover all vertices. 
However, these cycles may not be of distinct colours. 
First we show that for all $k \geq 3$, there are arbitrarily large $2$-edge-coloured complete $k$-graphs that cannot be partitioned into two monochromatic tight cycles of distinct colours.

\begin{prop}
\label{prop:extremal}
For all $k \geq 3$ and $m \geq k+1$, there exists a $2$-edge-colouring of $K_{k(m+1)+1}^{(k)}$ that does not admit a partition into two tight cycles of distinct colours.
\end{prop}

It is natural to ask whether we can cover almost all vertices of a $2$-edge-coloured complete $k$-graph with two vertex-disjoint monochromatic tight cycles of distinct colours. 
The case when $k=3$ is affirmed in~\cite{Bustamante2017}.
Here, we show that this is true when $k=4$. 

\begin{thm}
\label{thm:1}
For every $\eps > 0$, there exists an integer~$n_1$ such that, for all $n \geq n_1$, every $2$-edge-coloured complete $4$-graph on~$n$ vertices contains two vertex-disjoint monochromatic tight cycles of distinct colours covering all but at most~$\eps n$ of the vertices.
\end{thm}

When $k=5$, we prove a weaker result that four monochromatic tight cycles are sufficient to cover almost all vertices. 

\begin{thm}
\label{thm:2}
For every $\eps > 0$, there exists an integer~$n_1$ such that, for all $n \geq n_1$, every $2$-edge-coloured complete $5$-graph on~$n$ vertices contains four vertex-disjoint monochromatic tight cycles covering all but at most~$\eps n$ of the vertices.
\end{thm}

To prove \cref{thm:1,thm:2}, we use the \emph{connected matching method} that has often been credited to {\L}uczak~\cite{Luczak1999}.
We now present a sketch proof for~\cref{thm:1}.
Consider a $2$-edge-coloured complete $4$-graph~$K_n^{(4)}$ on~$n$ vertices.
We start by applying the Hypergraph Regularity Lemma to the 2-edge-coloured complete $4$-graph~$K_n^{(4)}$. 
More precisely the Regular Slice Lemma of Allen, B\"ottcher, Cooley and Mycroft~\cite{Allen2017}, see \cref{lem:regular slice}.
We obtain a 2-edge-coloured reduced graph~$\mathcal{R}$ that is almost complete.
A monochromatic matching in a $k$-graph is a set of vertex-disjoint edges of the same colour.
We say that it is tightly connected if, for any two edges~$f$ and~$f'$, there exists a sequence of edges $e_1,\dots,e_t$ of the same colour such that~$e_1 = f$,~$e_t=f'$ and $\s{e_i \cap e_{i+1}} = k-1$ for all $i \in [t-1]$.
Using~\cref{cor:matchings_to_cycles}, it suffices to find two vertex-disjoint monochromatic tightly connected matchings of distinct colours in the reduced graph $\mathcal{R}$.
The main challenge is to identify the `tightly connected components' (see \cref{section:preliminaries} for the formal definition) in which we will find the matchings. 
To do so, we introduce the concept of `blueprint', which is a $2$-edge-coloured $2$-graph with the same vertex set as~$\mathcal{R}$.
The key property is that connected components in the blueprint correspond to tightly connected components in~$\mathcal{R}$.

We conclude the introduction by outlining the structure of the paper.
In \cref{section:preliminaries}, we introduce some basic notation and definitions.
In \cref{section:extremal}, we prove \cref{prop:extremal}.
In \cref{section:regularity}, we introduce the statements about hypergraph regularity and prove the crucial \cref{cor:matchings_to_cycles} that allows us to reduce our problem of finding cycles in the complete graph to one about finding tightly connected matchings in the reduced graph.
In \cref{section:blueprints}, we give the definition of blueprint and setup some useful results. 
In \cref{section:matchings Kn4,section:matchings Kn5}, we prove \cref{thm:1,thm:2}, respectively. 
Finally, we make some concluding remarks in \cref{sec:concluding}.

\section{Preliminaries}
\label{section:preliminaries}
If we say that a statement holds for $0 < a \ll b \leq 1$, then we mean that there exists a non-decreasing function $f\colon (0,1] \rightarrow (0,1]$ such that the statement holds for all $a, b \in (0,1]$ with $a \leq f(b)$. Similar expressions with more variables are defined analogously. If~$1/n$ appears in one of these expressions, then we implicitly assume that~$n$ is a positive integer.

We often write $x_1\dots x_j$ for the set $\{x_1, \dots, x_j\}$. Moreover, for each positive integer~$n$, we let $[n] = \{1, \dots, n\}$.

Throughout this paper, any $2$-edge-colouring uses the colours red and blue. 
Let~$H$ be a $2$-edge-coloured $k$-graph. 
We denote by~$H^{\red}$ (and~$H^{\blue}$) the subgraph of~$H$ on~$V(H)$ induced by the red (and blue) edges of~$H$. 
Two edges~$f$ and~$f'$ in~$H$ are \emph{tightly connected} if there exists a sequence of edges $e_1,\dots,e_t$ such that~$e_1 = f$,~$e_t=f'$ and $\s{e_i \cap e_{i+1}} = k-1$ for all $i \in [t-1]$.
A subgraph~$H'$ of~$H$ is \emph{tightly connected} if every pair of edges in~$H'$ is tightly connected in~$H$. 
A maximal tightly connected subgraph of~$H$ is called a \emph{tight component} of~$H$.
Note that a tight component is a subgraph rather than a vertex subset as in the traditional graph case. 
A \emph{red tight component} and a \emph{red tightly connected matching} are a tight component and a tightly connected matching in~$H^{\red}$, respectively.
We define these terms similarly for \emph{blue}.

Let~$H$ be a $k$-graph and $S, W \subseteq V(H)$.
We denote by~$H-W$ the $k$-graph with $V(H-W) = V(H) \setminus W$ and $E(H-W) = \left\{e \in E(H) \colon e \cap W = \varnothing \right\}$. 
We call~$H-W$ the \emph{$k$-graph obtained from~$H$ by deleting~$W$}.
Further we let $H[W] = H - (V(H)\setminus W)$. 
Let~$F$ be a $k$-graph or a set of $k$-element sets. 
We denote by~$H-F$ the subgraph of~$H$ obtained by deleting the edges in~$F$. 
We define~$N_H(S,W)$ to be the set $\{e \in \binom{W}{k-|S|}\colon e\cup S \in H \}$ and we define~$d_H(S,W)$ to be its cardinality.
Further we write~$N_H(S)$ and~$d_H(S)$ for $N_H(S,V(H))$ and $d_H(S,V(H))$, respectively. 
If~$H$ is $2$-edge-coloured, then we write $N_H^{\red}(S,W)$, $d_H^{\red}(S,W)$, $N_H^{\blue}(S,W)$, $d_H^{\blue}(S,W)$ for $N_{H^{\red}}(S,W)$, $d_{H^{\red}}(S,W)$, $N_{H^{\blue}}(S,W)$, $d_{H^{\blue}}(S,W)$, respectively.
The \emph{link graph of~$H$ with respect to}~$S$, denoted by~$H_S$, is the $(k-\s{S})$-graph satisfying $V(H_S) = V(H)\setminus S$ and $E(H_S) = N_H(S)$.

For $j \in [k-1]$, the \emph{$j$-th shadow of~$H$}, denoted by~$\partial^j H$, is the \mbox{$(k-j)$}-graph with vertex set $V(\partial^j H) = V(H)$ and edge set \[E(\partial^j H) = \left\{e \in \binom{V(H)}{k-j}\colon e \subseteq f \text{ for some } f \in E(H) \right\}. \]
For the $1$-st shadow of~$H$, we also simply write $\partial H$ instead of $\partial^1H$.

For $\mu,\alpha > 0$, we say that a $k$-graph~$H$ on~$n$ vertices is $(\mu,\alpha)$\emph{-dense} if, for each $i \in [k-1]$, we have $d_H(S) \geq \mu \binom{n}{k-i}$ for all but at most~$\alpha\binom{n}{i}$ sets $S \in \binom{V(H)}{i}$ and $d_H(S) = 0$ for all other $S \in \binom{V(H)}{i}$. 

\begin{prop}
\label{prop:edges_in_dense}
Let $0 \leq \alpha , \mu \leq1$ and let~$H$ be a $(\mu, \alpha)$-dense $k$-graph on~$n$ vertices. Then $\s{H} \geq (\mu - \alpha)\binom{n}{k}$. Moreover, if $\mu > 1/2$, then~$H$ is tightly connected. 
\end{prop}
\begin{proof}
Note that 
\begin{align*}
\s{H} = \frac{1}{k} \sum_{S \in \binom{V(H)}{k-1}} d_H(S) \geq \frac{1}{k} (1-\alpha)\binom{n}{k-1} \mu n \geq (\mu - \alpha) \binom{n}{k}.
\end{align*}
Now suppose that $\mu > 1/2$. 
We show that~$H$ is tightly connected. Note that, for $S, S' \in \binom{V(H)}{k-1}$ with $d_H(S), d_H(S') > 0$, we have $d_H(S), d_H(S') \geq \mu n > n/2$ and thus 
\begin{align*}
N_H(S) \cap N_H(S') \neq \varnothing.
\end{align*}
Let $f = x_1 \dots x_k$ and $f' = y_1 \dots y_k$ be two edges of~$H$. Iteratively choose vertices $z_1, \dots, z_{k-1} \in V(H)$ such that 
\[
z_i \in N_H(z_1 \dots z_{i-1} x_{i+1} \dots x_k) \cap N_H(z_1 \dots z_{i-1} y_{i+1} \dots y_k)
\]
for all $i \in [k-1]$. 
It follows that~$f$ and~$f'$ are tightly connected. 
\end{proof}

The following proposition shows that any $k$-graph that has all but a small fraction of the possible edges contains a $(1- \varepsilon, \alpha)$-dense subgraph.
The proof was inspired by the proof of Lemma 8.8 in \cite{Han2017}. A different generalisation of this lemma can also be found as Lemma~$2.3$ in \cite{Lang2020}.

\begin{prop}
\label{prop:dense}
Let $1/n \ll \alpha \ll 1/k \leq 1/2$.
Let~$H$ be a $k$-graph on~$n$ vertices with $\s{H} \geq (1-\alpha)\binom{n}{k}$. Then there exists a subgraph~$H'$ of~$H$ such that $V(H') = V(H)$ and~$H'$ is $(1-2\alpha^{1/4k^2}, 2\alpha^{1/4k})$-dense. 
\end{prop}
\begin{proof}
We call a set $S \subseteq V(H)$ with $\s{S} \in [k-1]$ \emph{bad} if $d_H(S) < (1-\alpha^{1/2})\binom{n}{k-\s{S}}$. For $i \in [k-1]$, let~$\mathcal{B}_i$ be the set of all bad $i$-sets. For each $i \in [k-1]$, we have
\[
(1-\alpha)\binom{k}{i}\binom{n}{k} \leq \binom{k}{i} \s{H} = \sum_{S \in \binom{V(H)}{i}} d_H(S) \leq \binom{n}{i} \binom{n}{k-i} - \alpha^{1/2}\binom{n}{k-i}\s{\mathcal{B}_i}.
\]
This implies
\[
\s{\mathcal{B}_i} \leq \frac{1}{\alpha^{1/2}} \left( \binom{n}{i} - \frac{(1-\alpha)\binom{k}{i}\binom{n}{k}}{\binom{n}{k-i}} \right) \leq 2 \alpha^{1/2} \binom{n}{i}.
\]
Let $\beta = \alpha^{1/2k}$. For all $j \in \{k-1, k-2, \dots, 1\}$ in turn, we construct $\mathcal{A}_j \subseteq \binom{V(H)}{j}$ inductively as follows. We set $\mathcal{A}_{k-1} = \mathcal{B}_{k-1}$. Given $2 \leq j \leq k-1$ and $\mathcal{A}_j$, we define $\mathcal{A}_{j-1} \subseteq \binom{V(H)}{j-1}$ to be the set of all $X \in \binom{V(H)}{j-1}$ such that $X \in \mathcal{B}_{j-1}$ or $d_{\mathcal{A}_j}(X) \geq \beta^{1/2}n$. 
\begin{claim}
\label{claim:dense1}
For all $i \in [k-1]$, $\s{\mathcal{A}_i} \leq \beta^i \binom{n}{i}$. Moreover, if $1 \leq i < j \leq k-1$ and a set $S \in \binom{V(H)}{i}$ satisfies $d_{\mathcal{A}_j}(S) \geq \beta^{1/2(j-i)}\binom{n}{j-i}$, then $S \in \mathcal{A}_i$. 
\end{claim}
\begin{proofclaim}
We first prove the first part by induction on~$k-i$. For $i = k-1$, we have $\s{\mathcal{A}_{k-1}} = \s{\mathcal{B}_{k-1}} \leq 2 \alpha^{1/2}\binom{n}{k-1} \leq \beta^{k-1} \binom{n}{k-1}$.

Now suppose $2 \leq i \leq k-1$ and $\s{\mathcal{A}_i} \leq \beta^i\binom{n}{i}$.
By double counting tuples~$(X,w)$ with $X \in \mathcal{A}_{i-1} \setminus \mathcal{B}_{i-1}$ and $X \cup w \in \mathcal{A}_i$, we have $\br{\s{\mathcal{A}_{i-1}}- \s{\mathcal{B}_{i-1}}}\beta^{1/2}n \leq i \s{\mathcal{A}_i}$. Hence
\begin{align*}
\s{\mathcal{A}_{i-1}} &\leq \frac{i}{\beta^{1/2}n}\s{\mathcal{A}_i} + \s{\mathcal{B}_{i-1}} 
\leq \frac{i}{\beta^{1/2}n}\beta^i\binom{n}{i} + 2 \alpha^{1/2}\binom{n}{i-1} \\
&= \beta^{i-1/2} \binom{n-1}{i-1} + 2 \alpha^{1/2}\binom{n}{i-1} 
\leq \beta^{i-1} \binom{n}{i-1}.
\end{align*}
This proves the first part of the claim.

We now prove the second part of the claim. Fix $i \in [k-1]$. We proceed by induction on~$j-i$. For $j = i+ 1$, the statement holds by the definition of $\mathcal{A}_i$.
Now let $S \in \binom{V(H)}{i}$ and $j \geq i+2$ be such that $d_{\mathcal{A}_j} (S) \geq \beta^{1/2(j-i)} \binom{n}{j-i}$.
If $S \in \mathcal{B}_i$, then $S \in \mathcal{A}_i$. Recall that if $T \in \binom{V(H)}{j-1} \setminus \mathcal{A}_{j-1}$, then $d_{\mathcal{A}_j} (T) < \beta^{1/2}n$. We have 
\begin{align*}
\beta^{1/2(j-i)}\binom{n}{j-i} 
& \leq d_{\mathcal{A}_j}(S) \leq \sum_{\substack{ T \in \mathcal{A}_{j-1} \\ S \subseteq T}} d_{\mathcal{A}_{j}}(T) + \sum_{\substack{ T \in \binom{V(H)}{j-1} \setminus \mathcal{A}_{j-1} \\ S \subseteq T}} d_{\mathcal{A}_{j}}(T) \\ 
&\leq n d_{\mathcal{A}_{j-1}}(S) + \beta^{1/2}n d_{\binom{V(H)}{j-1}\setminus \mathcal{A}_{j-1}}(S) \\
&\leq n d_{\mathcal{A}_{j-1}}(S) + \beta^{1/2}n \binom{n}{j-i-1},
\end{align*}
and thus
\[
d_{\mathcal{A}_{j-1}}(S) \geq \beta^{1/2(j-i-1)}\binom{n}{j-i-1}.
\]
Hence by the induction hypothesis we have $S \in \mathcal{A}_i$.
\end{proofclaim}
For each $j \in [k-1]$, let~$F_j$ be the set of edges $e \in H$ for which there exists some $S \in \mathcal{A}_j$ with $S \subseteq e$. 
Let $F = \bigcup_{j \in [k-1]} F_j$ and let $H' = H - F$.
We will show that~$H'$ is the desired $k$-graph. 
For $i \in [k-1]$, let $\mathcal{S}_i$ be the set of all $S \in \binom{V(H)}{i}$ such that $d_F(S) \geq \beta^{1/2k}\binom{n}{k-i}$. 
\begin{claim}
For $i \in [k-1]$, $\s{\mathcal{S}_i} \leq \beta^{1/2}\binom{n}{i}$.
\end{claim}
\begin{proofclaim}
For $j \in [k-1]$, we have 
\[
\s{F_j} \leq \s{\mathcal{A}_j} \binom{n-j}{k-j} \overset{\text{\cref{claim:dense1}}}{\leq} \beta^j\binom{n}{j}\binom{n-j}{k-j} = \beta^j\binom{k}{j}\binom{n}{k}.
\]
Thus
\[
\s{F} \leq \sum_{j\in [k-1]} \s{F_j} \leq \sum_{j \in [k-1]} \beta^j \binom{k}{j} \binom{n}{k} \leq 2^k \beta \binom{n}{k}.
\]
Now, for $i \in [k-1]$, we have
\[
\frac{\s{\mathcal{S}_i}\beta^{1/2k}\binom{n}{k-i}}{\binom{k}{i}} \leq \s{F} \leq 2^k \beta \binom{n}{k}
\]
and thus $\s{\mathcal{S}_i} \leq \beta^{1/2}\binom{n}{i}$.
\end{proofclaim}

Consider $i \in [k-1]$. 
Note that $\s{\mathcal{S}_i \cup \mathcal{B}_i} \leq 2\alpha^{1/4k} \binom{n}{i}$. Now let $S \in \binom{V(H)}{i} \setminus \br{\mathcal{S}_i \cup \mathcal{B}_i}$. As $S \not\in \mathcal{B}_i$, we have $d_H(S) \geq (1-\alpha^{1/2}) \binom{n}{k-i}$.
As $S \not\in \mathcal{S}_i$, we have 
\begin{align*}
d_{H'}(S) &= d_H(S) - d_F(S) \geq d_H(S) - \beta^{1/2k}\binom{n}{k-i} \\ 
&\geq (1-\alpha^{1/2}- \beta^{1/2k}) \binom{n}{k-i} 
\geq (1 - 2\alpha^{1/4k^2})\binom{n}{k-i}.
\end{align*}
Consider $X \in \binom{V(H)}{i}$ with $d_{H'}(X) \neq 0$. 
We want to show that $d_{H'}(X) \geq (1- 2\alpha^{1/4k^2})\binom{n}{k-i}$. 
By the above, it suffices to show that $X \not\in \mathcal{B}_i \cup \mathcal{S}_i$.
Let $e \in H'$ with $X \subseteq e$. 
Since $e \not\in F_i$, we have $X \not\in \mathcal{A}_i$ and thus $X \not\in \mathcal{B}_i$. 
It remains for us to show that $X \not\in \mathcal{S}_i$. Assume the contrary that~$X$ is contained in more that $\beta^{1/2k} \binom{n}{k-i}$ edges of~$F$. 
Let $\mathcal{Y} = N_F(X)$, so $\s{\mathcal{Y}} \geq \beta^{1/2k} \binom{n}{k-i}$.
For each $Y \in \mathcal{Y}$, fix a set $A_Y \in \bigcup_{j \in [k-1]} \mathcal{A}_j$ such that $A_Y \subseteq X \cup Y$ and let $T_Y = X \cap A_Y$ and $S_Y = Y \setminus A_Y$. 
If $A_Y \subseteq X$, then $A_Y \subseteq e \in H'$, a contradiction.
Hence $A_Y \setminus X \neq \varnothing$ for all $Y \in \mathcal{Y}$. Thus, for $Y \in \mathcal{Y}$, we have $\s{T_Y} \leq \s{A_Y} -1 \leq k-2$.
By an averaging argument, there exist $t \in \{0, 1, \dots, k-2\}$, $T \in \binom{X}{t}$, $a \in [k-1]$, $S \in \binom{V(H)}{k-i-a+t}$ and $\tY \subseteq \mathcal{Y}$ such that, for all $Y \in \tY$, we have $T_Y = T$, $\s{A_Y} = a$, $S_Y = S$ and 
\[
\s{\tY} \geq \frac{\s{\mathcal{Y}}}{2^i (k-1) \binom{n}{k-i-a+t}} \geq \beta^{1/2(k-1)}\binom{n}{a-t}.
\]
Since $Y \setminus A_Y = S_Y = S$ and $\s{A_Y} = a$ for all $Y \in \tY$, the~$A_Y$ are distinct for all $Y \in \tY$. 
Recall that $T \subseteq A_Y \in \mathcal{A}_a$ for each $Y \in \tY$. 
If $T = \varnothing$, then $t=0$ and so $\s{\mathcal{A}_a} \ge \s{\tY} > \beta^a \binom{n}{a}$ contradicting \cref{claim:dense1}.
If $T \ne \varnothing$, then we have $d_{\mathcal{A}_a}(T) \ge \s{\tY} \ge \beta^{1/2(k-1)}\binom{n}{a-t}$. 
\cref{claim:dense1} implies that $T \in \mathcal{A}_t$.
Since $T \subseteq X \subseteq e$, we have $e \in F_t$ contradicting the fact that $e \in H' = H - \bigcup_{j \in [k-1]}F_j
$. 
\end{proof}

\section{Extremal example}
\label{section:extremal}

In this section, we prove \cref{prop:extremal}, that is, we prove that, for $k \geq 3$, there exist arbitrarily large $2$-edge-coloured complete $k$-graphs that do not admit a partition into two tight cycles of distinct colours.

A \emph{$k$-uniform tight path} is a $k$-graph obtained by deleting a vertex from a tight cycle.
First we need the following proposition.
\begin{prop}
\label{prop:path}
Let $k \geq 3$, let~$P$ and~$C$ be a $k$-uniform tight path and tight cycle, respectively. We have the following.
\begin{enumerate}[label = \upshape (\roman*)]
    \item \label{extrprop1} If~$X$ and~$Y$ partition~$V(P)$ such that $\s{e \cap Y} \geq 2$ for all $e \in P$, then $2(\s{X} - (k-1)) \leq (k-2) \s{Y}.$
   
   \item \label{extrprop2} If~$X$ and~$Y$ partition~$V(C)$ such that $\s{e \cap Y} \geq 2$ for all $e \in C$, then $2\s{X} \leq (k-2) \s{Y}.$
\end{enumerate}
\end{prop}
\begin{proof}
We first prove \ref{extrprop1}.
Let~$M$ be a matching of maximum size in~$P$. Since each edge of~$P$ contains at least~$2$ vertices of~$Y$,
\[
\s{X} \leq \s{X \cap V(M)} + \s{V(P) \setminus V(M)} \leq (k-2)\s{M} + k-1 \leq \frac{(k-2)\s{Y}}{2} + k-1.
\]

Now we prove \ref{extrprop2}. Since $\s{e \cap Y} \geq 2$ and $\s{e \cap X} \leq k-2$ for each edge $e \in C$, we have
\[
\s{X} = \frac{1}{k} \sum_{e \in C} \s{e \cap X} = \frac{1}{k} \sum_{e \in C} \frac{\s{e \cap X}}{\s{e \cap Y}} \s{e \cap Y} \leq \frac{1}{k} \sum_{e \in C} \frac{k-2}{2} \s{e \cap Y} = \frac{k-2}{2}\s{Y}. 
\] 
\end{proof}

We are now ready to give our extremal example.
Note that the case~$k=3$ of the extremal example is already given in \cite{Garbe2019}.
 Recall that, in a $k$-graph, we consider a single edge and any set of fewer than~$k$ vertices to be degenerate cycles.

\begin{proof}[Proof of \cref{prop:extremal}]
Let $k\geq 3$, $m \geq k+1$ and $n = k(m+1)+1$. Let~$X$,~$Y$ and~$\{z\}$ be three disjoint vertex sets of~$K_n^{(k)}$ of sizes~$(k-1)m+k-2$,~$m+2$ and 1, respectively. 
We colour an edge~$e$ in~$K_n^{(k)}$ red if $z \in e$ and $\s{e \cap Y} \geq 2$ or $z \not\in e$ and $\s{e \cap Y} = 1$. Otherwise we colour it blue.
Note that $K_n^{(k)} - z$ has the following~$3$ monochromatic tight components: 
\begin{align*}
B_1 = \binom{X}{k},\ 
B_2 = \left\{e \in \binom{X \cup Y}{k}\colon \s{e \cap Y} \geq 2\right\},\ 
R = \left\{e \in \binom{X\cup Y}{k} \colon \s{e\cap Y} =1 \right\}.
\end{align*}
Note that~$B_1$ and~$B_2$ are blue and~$R$ is red.
Suppose for a contradiction that~$K_n^{(k)}$ can be partitioned into a red tight cycle~$C_R$ and a blue tight cycle~$C_B$.

First assume $z \in V(C_R)$. Since all the red edges containing~$z$ are in a red tight component disjoint from~$R$, we have $\s{V(C_R)} \leq k.$ Hence $\s{V(C_B)} = n - \s{V(C_R)} \geq n-k \geq km > k$ and $\s{V(C_B) \cap Y} = \s{Y \setminus V(C_R)} \geq m+2-(k-1) \geq 1$. 
So~$C_B$ is not degenerate and $C_B \subseteq B_2$.
Any edge~$e \in C_B$ contains at least~$2$ vertices in~$Y$. 
By \cref{prop:path}\ref{extrprop2}, $2 \s{V(C_B) \cap X} \leq (k-2)\s{V(C_B) \cap Y}$. It follows that 
\begin{align*}
2(k-1)m -2 &= 2(\s{X} - (k-1)) \leq 2 \s{V(C_B)\cap X} \\ &\leq (k-2)\s{V(C_B) \cap Y} \leq (k-2)\s{Y} = (k-2)(m+2).
\end{align*}
This implies that $m \leq 2$, a contradiction.

Hence, we may assume that $z \in V(C_B)$. This implies that $C_R \subseteq R$ or $\s{V(C_R)} \leq k-1$.
Let $x_R = \s{V(C_R) \cap X}$, $y_R = \s{V(C_R) \cap Y}$, $x_B = \s{V(C_B) \cap X}$ and $y_B = \s{V(C_B) \cap Y}$. Let~$P_B$ be the tight path $C_B -z$. Clearly $\s{V(P_B) \cap X} = x_B$ and $\s{V(P_B) \cap Y} = y_B$.
Since $C_R \subseteq R$ or $\s{V(C_R)} \leq k-1$,
\begin{align}
\label{y_R_bound}
y_R \leq \max \left\{ \left\lfloor \frac{\s{X}}{k-1}\right\rfloor, k-1 \right\} = m < \s{Y}.
\end{align}
Hence, $V(P_B) \cap Y \neq \varnothing$ and $\s{V(P_B)} \geq (n-1) - km \geq k$.
We must have $P_B \subseteq B_2$. By \cref{prop:path}\ref{extrprop1}, we have that 
\begin{align}
\label{eqxb}
2(x_B - (k-1)) \leq (k-2)y_B.
\end{align}
Thus
\begin{align*}
\s{V(P_B)} 
&= x_B + y_B \leq \frac{k}{2}y_B + k-1 
\leq \frac{k}{2}\s{Y} + k-1 
= \frac{k}{2}(m+2) + k-1 \\
&\leq mk = n-1 - k.
\end{align*}
This implies that $\s{V(C_R)} \geq k$. Hence $C_R \subseteq R$ and thus 
\begin{align}
\label{eqxr}
    x_R = (k-1)y_R.
\end{align}
Since $x_R + x_B = \s{X} = (k-1)m+k-2$ and $y_R + y_B = \s{Y} = m+2$, \cref{eqxb} implies
\begin{align*}
    (k-2)(m+2 - y_R) &\geq 2(\s{X} - x_R -(k-1)) \\
    &= 2((k-1)m+k-2-(k-1)y_R-(k-1)),
\end{align*}
which implies $y_R \geq m-1$.
If $y_R = m-1$, then \eqref{eqxr} implies that $x_R = (k-1)(m-1)$ and thus $x_B = 2k-3$ and~$y_B=3$. 
Let $P_B = v_1\dots v_{2k}$.
Either the edge $v_1 \dots v_k$ or the edge $v_{k+1}\dots v_{2k}$ contains at most one vertex of~$Y$, a contradiction to $P_B \subseteq B_2$.
Thus we may assume $y_R \geq m$ and since $y_R \leq m$ by \cref{y_R_bound}, we have $y_R = m$. 
By \cref{eqxr}, we have $x_R = (k-1)m$ and thus $x_B =k-2$ and $y_B = 2$.
Hence,~$C_B$ is a copy of $K_{k+1}^{(k)}$ that has a blue edge containing~$z$ and at least two vertices of~$Y$, a contradiction.
\end{proof}

\section{Hypergraph regularity}
 \label{section:regularity}

In this section, we give the formulation of hypergraph regularity that we use, following closely the presentation of Allen, B\"ottcher, Cooley and Mycroft \cite{Allen2017}.
A \emph{hypergraph} $\mathcal{H}$ is an ordered pair $(V(\mathcal{H}), E(\mathcal{H}))$, where $E(\mathcal{H}) \subseteq 2^{V(\mathcal{H})}$. 
Again, we identify the hypergraph $\mathcal{H}$ with its edge set $E(\mathcal{H})$. 
A subgraph $\mathcal{H}'$ of $\mathcal{H}$ is a hypergraph with $V(\mathcal{H}') \subseteq V(\mathcal{H})$ and $E(\mathcal{H}') \subseteq E(\mathcal{H})$. 
It is \emph{spanning} if $V(\mathcal{H}') = V(\mathcal{H})$.
For $U \subseteq V(\mathcal{H})$, we define $\mathcal{H}[U]$ to be the subgraph of $\mathcal{H}$ with $V(\mathcal{H}[U]) = U$ and $E(\mathcal{H}[U]) = \{ e \in E(\mathcal{H}) \colon e \subseteq U\}$.
We call $\mathcal{H}$ a \emph{complex} if $\mathcal{H}$ is down-closed, that is if $e \in \mathcal{H}$ and $f \subseteq e$, then $f \in \mathcal{H}$.
A \emph{$k$-complex} is a complex with only edges of size at most~$k$. 
We denote by $\mathcal{H}^{(i)}$ the spanning subgraph of~$\mathcal{H}$ containing only the edges of size~$i$.
Let $\mathcal{P}$ be a partition of $V(\mathcal{H})$ into parts $V_1, \dots, V_s$. Then we say that a set $S \subseteq V(\mathcal{H})$ is \emph{$\mathcal{P}$-partite} if $\s{S \cap V_i} \leq 1$ for all $i \in [s]$. 
For $\mathcal{P}' = \{V_{i_1}, \dots, V_{i_r}\} \subseteq \mathcal{P}$, we define the subgraph of $\mathcal{H}$ induced by $\mathcal{P}'$, denoted by $\mathcal{H}[\mathcal{P'}]$ or $\mathcal{H}[V_{i_1}, \dots, V_{i_r}]$, to be the subgraph of $\mathcal{H}[\bigcup \mathcal{P}']$ containing only the edges that are $\mathcal{P}'$-partite. 
The hypergraph $\mathcal{H}$ is $\mathcal{P}$-partite if all of its edges are $\mathcal{P}$-partite. In this case we call the parts of $\mathcal{P}$ the \emph{vertex classes} of $\mathcal{H}$. We say that $\mathcal{H}$ is \emph{$s$-partite} if it is $\mathcal{P}$-partite for some partition $\mathcal{P}$ of $V(\mathcal{H})$ into~$s$ parts.
Let $\mathcal{H}$ be a $\mathcal{P}$-partite hypergraph. 
If~$X$ is a $k$-set of vertex classes of~$\mathcal{H}$, then we write~$\mathcal{H}_X$ for the $k$-partite subgraph of $\mathcal{H}^{(k)}$ induced by~$\bigcup X$, whose vertex classes are the elements of~$X$. 
Moreover, we denote by $\mathcal{H}_{X^<}$ the $k$-partite hypergraph with $V(\mathcal{H}_{X^<}) = \bigcup X$ and $E(\mathcal{H}_{X^<}) = \bigcup_{X'\subsetneq X} \mathcal{H}_{X'}$.
In particular, if $\mathcal{H}$ is a complex, then $\mathcal{H}_{X^<}$ is a $(k-1)$-complex because~$X$ is a set of size~$k$.

Let $i \geq 2$, and let $\mathcal{P}_i$ be a partition of a vertex set~$V$ into~$i$ parts. Let~$H_i$ and~$H_{i-1}$ be a $\mathcal{P}_i$-partite $i$-graph and a $\mathcal{P}_i$-partite $(i-1)$-graph on a common vertex set~$V$, respectively. We say that a $\mathcal{P}_i$-partite $i$-set in~$V$ is \emph{supported on}~$H_{i-1}$ if it induces a copy of the complete $(i-1)$-graph $K_i^{(i-1)}$ on~$i$ vertices in~$H_{i-1}$. We denote by $K_i(H_{i-1})$ the $\mathcal{P}_i$-partite $i$-graph on~$V$ whose edges are all $\mathcal{P}_i$-partite $i$-sets contained in~$V$ which are supported on~$H_{i-1}$. Now we define the \emph{density of~$H_i$ with respect to~$H_{i-1}$} to be
\[
d(H_i \mid H_{i-1}) = \frac{\s{K_i(H_{i-1}) \cap H_i}}{\s{K_i(H_{i-1})}}
\]
if $\s{K_i(H_{i-1})} > 0$ and $d(H_i \mid H_{i-1}) = 0$ if $\s{K_i(H_{i-1})} = 0$.
So $d(H_i \mid H_{i-1})$ is the proportion of $\mathcal{P}_i$-partite copies of~$K_i^{i-1}$ in~$H_{i-1}$ which are also edges of~$H_i$. More generally, if $\mathbf{Q} = (Q_1, Q_2, \dots, Q_r)$ is a collection of~$r$ (not necessarily disjoint) subgraphs of~$H_{i-1}$, we define $K_i(\mathbf{Q}) = \bigcup_{j=1}^r K_i(Q_j)$ and 
\[
d(H_i \mid \mathbf{Q}) = \frac{\s{K_i(\mathbf{Q}) \cap H_i}}{\s{K_i(\mathbf{Q})}}
\]
if $\s{K_i(\mathbf{Q})} > 0$ and $d(H_i \mid \mathbf{Q}) = 0$ if $\s{K_i(\mathbf{Q})} = 0$.
We say that~$H_i$ is \emph{$(d_i, \eps, r)$-regular with respect to~$H_{i-1}$}, if we have $d(H_i \mid \mathbf{Q}) = d_i \pm \eps$ for every $r$-set~$\mathbf{Q}$ of subgraphs of~$H_{i-1}$ with $\s{K_i(\mathbf{Q})} > \eps \s{K_i(H_{i-1})}$.
We say that~$H_i$ is \emph{$(\eps, r)$-regular with respect to~$H_{i-1}$} if there exists some~$d_i$ for which~$H_i$ is $(d_i, \eps, r)$-regular with respect to~$H_{i-1}$. 
Finally, given an $i$-graph~$G$ whose vertex set contains that of~$H_{i-1}$, we say that~$G$ is \emph{$(d_i, \eps, r)$-regular with respect to~$H_{i-1}$} if the $i$-partite subgraph of~$G$ induced by the vertex classes of~$H_{i-1}$ is $(d_i, \eps, r)$-regular with respect to~$H_{i-1}$. 
We refer to the density of this $i$-partite subgraph of~$G$ with respect to~$H_{i-1}$ as the \emph{relative density of~$G$ with respect to~$H_{i-1}$}.

Now let $s \geq k \geq 3$ and let $\mathcal{H}$ be an $s$-partite $k$-complex on vertex classes $V_1, \dots, V_s$. For any set $A \subseteq [s]$, we write~$V_A$ for $\bigcup_{i \in A} V_i$. Note that, if $e \in \mathcal{H}^{(i)}$ for some $2 \leq i \leq k$, then the vertices of~$e$ induce a copy of~$K_i^{i-1}$ in $\mathcal{H}^{(i-1)}$. Therefore, for any set $A \in \binom{[s]}{i}$, the density $d(\mathcal{H}^{(i)}[V_A] \mid \mathcal{H}^{(i-1)}[V_A])$ is the proportion of `possible edges' of $\mathcal{H}^{(i)}[V_A]$, which are indeed edges. We say that $\mathcal{H}$ is \emph{$(d_k, \dots, d_2, \eps_k, \eps, r)$-regular} if
\begin{enumerate}[label=(\alph*)]
    \item for any $2 \leq i \leq k-1$ and any $A \in \binom{[s]}{i}$, the induced subgraph $\mathcal{H}^{(i)}[V_A]$ is $(d_i, \eps, 1)$-regular with respect to $\mathcal{H}^{(i-1)}[V_A]$, and 
    \item for any $A \in \binom{[s]}{k}$, the induced subgraph $\mathcal{H}^{(k)}[V_A]$ is $(d_k, \eps_k, r)$-regular with respect to $\mathcal{H}^{(k-1)}[V_A]$.
\end{enumerate}
For $\mathbf{d} = (d_k, \dots, d_2)$, we write $(\mathbf{d}, \eps_k, \eps, r)$-regular to mean $(d_k, \dots, d_2, \eps_k, \eps, r)$-regular.
We say that a $(k-1)$-complex $\mathcal{J}$ is \emph{$(t_0, t_1, \eps)$-equitable} if it has the following properties.
\begin{enumerate}[label = (\alph*)]
    \item $\mathcal{J}$ is $\mathcal{P}$-partite for some $\mathcal{P}$ which partitions $V(\mathcal{J})$ into~$t$ parts, where $t_0 \leq t \leq t_1$, of equal size. We refer to $\mathcal{P}$ as the \emph{ground partition} of $\mathcal{J}$, and to the parts of $\mathcal{P}$ as the \emph{clusters} of $\mathcal{J}$.
    \item There exists a \emph{density vector} $\mathbf{d} = (d_{k-1}, \dots, d_2)$ such that, for each $2 \leq i \leq k-1$, we have $d_i \geq 1/t_1$ and $1/d_i \in \mathbb{N}$, and $\mathcal{J}$ is $(\mathbf{d}, \eps, \eps, 1)$-regular.
\end{enumerate}
For any $k$-set~$X$ of clusters of $\mathcal{J}$, we denote by $\hat{\mathcal{J}}_X$ the $k$-partite $(k-1)$-graph $(\mathcal{J}_{X^<})^{(k-1)}$ and call $\hat{\mathcal{J}}_X$ a \emph{polyad}. Given a $(t_0, t_1, \eps)$-equitable $(k-1)$-complex $\mathcal{J}$ and a $k$-graph~$G$ on $V(\mathcal{J})$, we say that~$G$ is \emph{$(\eps_k, r)$-regular with respect to a $k$-set~$X$ of clusters of $\mathcal{J}$} if there exists some~$d$ such that~$G$ is $(d, \eps_k, r)$-regular with respect to the polyad $\hat{\mathcal{J}}_X$.
Moreover, we write $d_{G, \mathcal{J}}^*(X)$ for the relative density of~$G$ with respect to $\hat{\mathcal{J}}_X$; we may drop either subscript if it is clear from context. 

We can now give the crucial definition of a regular slice.
\begin{defn}[Regular slice]
Given $\eps, \eps_k > 0, r, t_0, t_1 \in \mathbb{N}$, a graph~$G$ and a $(k-1)$-complex $\mathcal{J}$ on~$V(G)$, we call $\mathcal{J}$ a \emph{$(t_0, t_1, \eps, \eps_k,r)$-regular slice} for~$G$ if $\mathcal{J}$ is $(t_0, t_1, \eps)$-equitable and~$G$ is $(\eps_k, r)$-regular with respect to all but at most $\eps_k \binom{t}{k}$ of the $k$-sets of clusters of~$\mathcal{J}$, where~$t$ is the number of clusters of $\mathcal{J}$.
\end{defn}

If we specify the density vector~$\mathbf{d}$ and the number of clusters~$t$ of an equitable complex or a regular slice, then it is not necessary to specify~$t_0$ and~$t_1$ (since the only role of these is to bound~$\mathbf{d}$ and~$t$). In this situation we write that $\mathcal{J}$ is $(\cdot, \cdot, \eps)$-equitable, or is a $(\cdot, \cdot, \eps, \eps_k, r)$-regular slice for~$G$.

Given a regular slice $\mathcal{J}$ for a $k$-graph~$G$, we define the $d$-reduced $k$-graph $\mathcal{R}_d^{\mathcal{J}}(G)$ as follows.
\begin{defn}[The $d$-reduced $k$-graph]
Let $k \geq 3$.
Let~$G$ be a $k$-graph and let $\mathcal{J}$ be a $(t_0, t_1, \eps, \eps_k, r)$-regular slice for~$G$. Then, for~$d >0$, we define the \emph{$d$-reduced $k$-graph $\mathcal{R}_d^{\mathcal{J}}(G)$} to be the $k$-graph whose vertices are the clusters of $\mathcal{J}$ and whose edges are all $k$-sets~$X$ of clusters of $\mathcal{J}$ such that~$G$ is $(\eps_k, r)$-regular with respect to~$X$ and $d^*(X) \geq d$.
\end{defn}

We now state the version of the Regular Slice Lemma that we need, which is a special case of~\cite[Lemma 10]{Allen2017}.

\begin{lem}[Regular Slice Lemma {\cite[Lemma 10]{Allen2017}}]
\label{lem:regular slice}
Let $k \geq 3$. For all positive integers~$t_0$ and~$s$, positive~$\eps_k$ and all functions $r\colon \mathbb{N} \rightarrow \mathbb{N}$ and $\eps \colon \mathbb{N} \rightarrow (0,1]$, there are integers~$t_1$ and~$n_0$ such that the following holds for all $n \geq n_0$ which are divisible by~$t_1!$. 
Let~$K$ be a $2$-edge-coloured complete $k$-graph on~$n$ vertices.
Then there exists a $(k-1)$-complex~$\mathcal{J}$ on~$V(K)$ which is a $(t_0, t_1, \eps(t_1), \eps_k, r(t_1))$-regular slice for both~$K^{\red}$ and~$K^{\blue}$.
\end{lem}

Given a $2$-edge-coloured complete $k$-graph~$H$ we want to apply the Regular Slice Lemma to~$H^{\red}$ and~$H^{\blue}$. The following lemma shows that in this setting the union of the corresponding reduced graphs $\mathcal{R}_d^\mathcal{J}(H^{\red}) \cup \mathcal{R}_d^\mathcal{J}(H^{\blue})$ is almost complete.

\begin{lem}[{\cite[Lemma 8.5]{Garbe2019}}]
\label{lem:reduced graph edge count}
Let $k \geq 3$.
Let~$K$ be a $2$-edge-coloured complete $k$-graph and let $\mathcal{J}$ be a $(\cdot,\cdot,\eps, \eps_k,r)$-regular slice for both~$K^{\red}$ and~$K^{\blue}$. Let~$t$ be the number of clusters of $\mathcal{J}$. Then, provided that $d \leq 1/2$, 
we have $\s{\mathcal{R}_d^{\mathcal{J}}(K^{\red}) \cup \mathcal{R}_d^{\mathcal{J}}(K^{\blue})} \geq (1 - 2\eps_k)\binom{t}{k}$.
\end{lem}
\begin{proof}
Since $\mathcal{J}$ is a $(\cdot,\cdot,\eps, \eps_k,r)$-regular slice for both~$K^{\red}$ and~$K^{\blue}$ there are at least $(1-2 \eps_k)\binom{t}{k}$ $k$-sets~$X$ of clusters of $\mathcal{J}$ such that both~$K^{\red}$ and~$K^{\blue}$ are $(\eps_k, r)$-regular with respect to~$X$. Let~$X$ be such a $k$-set. Since~$K^{\red}$ and~$K^{\blue}$ are complements of each other, we have $d_{K^{\red}}^*(X) + d_{K^{\blue}}^*(X) = 1$. Hence $d_{K^{\red}}^*(X) \geq 1/2$ or $d_{K^{\blue}}^*(X) \geq 1/2$ and thus, since $d \leq 1/2$, we have $X \in \mathcal{R}_d^{\mathcal{J}}(K^{\red}) \cup \mathcal{R}_d^{\mathcal{J}}(K^{\blue})$.
\end{proof}

Let~$H$ be a $k$-graph. 
A \emph{fractional matching} in~$H$ is a function $\omega : E(H) \rightarrow [0, 1]$ such that for all $v \in V(H)$, $\omega(v) \coloneqq \sum_{e \in H: v \in e} \omega(e) \le 1$. 
The \emph{weight} of the fractional matching is defined to be $\sum_{e \in H} \omega(e)$. 
A fractional matching is \emph{tightly connected} if the subgraph induced by the edges~$e$ with $\omega (e)>0$ is tightly connected in~$H$. 
The following result from~{\cite{Allen2017}} converts a tightly connected fractional matching in the reduced graph into a tight cycle in the original graph. 

\begin{lem}[{\cite[Lemma 13]{Allen2017}}]
\label{lem:cycle}
Let~$k,r,n_0,t$ be positive integers, and let $\psi, \eps, \eps_k, d_k, \dots, d_2$ be positive constants such that $1/d_i \in \mathbb{N}$ for each $2 \leq i \leq k-1$, and such that $1/n_0 \ll 1/t$,
\[
\frac{1}{n_0} \ll \frac{1}{r}, \eps \ll \eps_k, d_{k-1}, \dots, d_2 \quad \text{and} \quad \eps_k \ll \psi, d_k, \frac{1}{k}.
\]
Then the following holds for all integers $n \geq n_0$. Let~$G$ be a $k$-graph on~$n$ vertices, and $\mathcal{J}$ be a $(\cdot, \cdot, \eps, \eps_k, r)$-regular slice for~$G$ with~$t$ clusters and density vector $(d_{k-1}, \dots, d_2)$. Suppose that $\mathcal{R}_{d_k}^{\mathcal{J}}(G)$ contains a tightly connected fractional matching with weight~$\mu$.
Then~$G$ contains a tight cycle of length~$\ell$ for every $\ell \leq (1-\psi)k\mu n/t$ that is divisible by~$k$.
\end{lem}

We use the following fact, lemma and proposition to prove \cref{lem:main} which is a stronger version of \cref{lem:cycle} that allows us to control the location of the tight cycle.

\begin{fact}[{\cite[Fact 7]{Allen2017}}]
\label{fact}
Suppose that $1/m_0 \ll \eps \ll 1/t_1, 1/t_0, \beta, 1/k \leq 1/3$ and that $\mathcal{J}$ is a $(t_0, t_1, \eps)$-equitable $(k-1)$-complex with density vector $(d_{k-1}, \dots, d_2)$ whose clusters each have size $m \geq m_0$. Let~$X$ be a set of~$k$ clusters of $\mathcal{J}$. Then
\[
\s{K_k((\mathcal{J}_{X^<})^{(k-1)})} = (1 \pm \beta)m^k \prod_{i=2}^{k-1} d_i^{\binom{k}{i}}.
\]
\end{fact}

\begin{lem}[Regular Restriction Lemma {\cite[Lemma 28]{Allen2017}}]
\label{regrestriction}
Suppose integers~$k,m$ and reals $\alpha, \eps, \eps_k, d_k, \dots, d_2 >0$ are such that 
\[
\frac{1}{m} \ll \eps \ll \eps_k, d_{k-1}, \dots, d_2 \quad \text{and} \quad \eps_k \ll \alpha, \frac{1}{k}.
\]
For any $r,s \in \mathbb{N}$ and~$d_k > 0$, set $\mathbf{d} = (d_k, \dots, d_2)$, and let $\mathcal{G}$ be an $s$-partite $k$-complex whose vertex classes $V_1, \dots, V_s$ each have size~$m$ and which is $(\mathbf{d}, \eps_k, \eps, r)$-regular. 
Choose any $V_i' \subseteq V_i$ with $\s{V_i'} \geq \alpha m$ for each $i \in [s]$. 
Then the induced subcomplex $\mathcal{G}[V_1'\cup \dots \cup V_s']$ is $(\mathbf{d}, \sqrt{\eps_k}, \sqrt{\eps}, r)$-regular.
\end{lem}

The following proposition shows that a refinement of a regular slice is also a regular slice. 

\begin{prop}
\label{prop:split}
Let $1/m \ll \eps \ll 1/N, 1/t_0, 1/t_1, 1/k \leq 1/3$.
Let $\mathcal{J}$ be a $(t_0,t_1, \eps)$-equitable $(k-1)$-complex with density vector $(d_{k-1}, \dots, d_2)$ and clusters $V_1, \dots, V_t$ each of size~$m$.
Let $V_{i,1}, \dots, V_{i,N}$ be an equipartition of~$V_i$ for each $i \in [t]$. Then there exists a $(N t_0, N t_1, \sqrt{\eps})$-equitable $(k-1)$-complex~$\J$ with density vector $(d_{k-1}, \dots, d_2)$, ground partition $\{V_{i,j} \colon i \in [t], j \in [N]\}$ and $\J[V_1, \dots, V_t] = \mathcal{J}$.
\end{prop}
\begin{proof}
We construct~$\J$ from $\mathcal{J}$ as follows.
Let the ground partition of~$\J$ be $\{V_{i,j} \colon i \in [t], j \in [N]\}$. Starting with the edges of $\mathcal{J}$ we iteratively add additional edges at random as follows. For each $2 \leq i \leq k-1$, beginning with~$i = 2$, we add each $i$-edge that contains two vertices that are in vertex classes with the same first index and is supported on the $(i-1)$-edges independently with probability~$d_i$.

We now show that with high probability~$\J$ is the desired $(k-1)$-complex.
Note that it suffices to show that with high probability~$\J$ is $(\mathbf{d},\sqrt{\eps},\sqrt{\eps},1)$-regular.

Let $\J^{\leq i} = \bigcup_{j \in [i]} \J^{(j)}$ and $\mathbf{d}^{\leq i} = (d_i, \dots, d_2)$.
For $i \in [k-1]$, let~$B_i$ be the event that $\J^{\leq i}$ is not $(\mathbf{d}^{\leq i}, \sqrt{\eps}, \sqrt{\eps}, 1)$-regular. Note that $B_1 = \varnothing$.
Consider $2 \leq i \leq k-1$ and $A \in \binom{[t] \times [N]}{i}$. Let~$B_{i,A}$ be the event that $\J^{(i)}[V_A]$ is not $(d_i, \sqrt{\eps}, 1)$-regular with respect to $\J^{(i-1)}[V_A]$.
\begin{claim}
\label{claim:probability_bound}
For $i \in [k-1]$ and $A \in \binom{[t] \times [N]}{i}$, we have
$\mathbb{P}\left[B_{i,A} \mid \overline{B_{i-1}} \,\right] = e^{-\Omega(m^i)}$ as $m \rightarrow \infty$.
\end{claim}
\begin{proofclaim}
Assume $\overline{B_{i-1}}$ holds.
Let $A = \{(r_j,s_j)\colon j \in [i]\}$. Define $\widetilde{A} = \{r_j \colon j \in [i]\}$. If the~$r_j$ are distinct, then the claim holds by \cref{regrestriction} with $\mathcal{G} = \mathcal{J}[V_{\widetilde{A}}]$ and $\alpha = 1/N$. 
If not all the~$r_j$ are distinct, then $\s{K_i(\J^{(i-1)}[V_A])} \geq \frac{1}{2}\br{\prod_{j=2}^{i-1}d_j^{\binom{i}{j}}}(m/N)^i$, by \cref{fact}. Thus for each subgraph~$Q$ of $\J^{(i-1)}[V_A]$ such that $\s{K_i(Q)} > \sqrt{\eps}\s{K_i(\J^{(i-1)}[V_A])}$, a Chernoff bound implies that
\begin{align*}
    &\mathbb{P}\left.\left[d(\J^{(i)}[V_A] \mid Q) \neq d_i \pm \sqrt{\eps} \right\vert \overline{B_{i-1}} \,\right]
    \\ = &\mathbb{P}\left.\left[ \s{\s{\J^{(i)}[V_A] \cap K_i(Q)} - d_i\s{K_i(Q)}} > \frac{\sqrt{\eps}}{d_i} d_i \s{K_i(Q)} \right\vert \overline{B_{i-1}} \,\right]
    \\ \leq &2\exp\br{-\frac{1}{3}\br{\frac{\sqrt{\eps}}{d_i}}^2d_i \s{K_i(Q)}}
     \leq 2\exp\br{-\frac{1}{3}\frac{\eps^{3/2}}{d_i}\s{K_i(\J^{(i-1)}[V_A]}}
    \\ \leq &2\exp\br{-\frac{1}{6}\frac{\eps^{3/2}}{d_i}\br{\prod_{j=2}^{i-1}d_j^{\binom{i}{j}}}\br{\frac{m}{N}}^i}
     \leq e^{-\Omega(m^i)}.
\end{align*}
Since there are at most $2^{(im)^{i-1}}$ choices for~$Q$, the claim follows by a union bound.
\end{proofclaim}
Note that if~$\J$ is not $(\mathbf{d}, \sqrt{\eps},\sqrt{\eps},1)$-regular, then there exists some $i \in [k-1]$ and $A \in \binom{[t] \times [N]}{i}$ such that~$B_{i,A}$ holds. Further by choosing~$i$ minimal we can ensure that $\overline{B_{i-1}}$ holds.
Thus, by a union bound and \cref{claim:probability_bound}, we have
\begin{align*}
    \mathbb{P}\left[\J \text{ is not $(\mathbf{d}, \sqrt{\eps},\sqrt{\eps},1)$-regular}\right] 
    \leq &\sum_{i=1}^{k-1}\sum_{A \in \binom{[t] \times [N]}{i}} \mathbb{P}\left[B_{i,A} \cap \overline{B_{i-1}} \, \right]
    \\ \leq &\sum_{i=1}^{k-1}\sum_{A \in \binom{[t] \times [N]}{i}} \mathbb{P}\left[B_{i,A}\mid \overline{B_{i-1}} \,\right] = o(1).
\end{align*}
\end{proof}

The following lemma is a strengthening of \cref{lem:cycle}. We believe the constant~$\beta$ and the corresponding condition could be removed if one were to go through the proof of \cref{lem:cycle} to prove a stronger result.

\begin{lem}
\label{lem:main}
Let $1/n \ll 1/r, \eps \ll \eps_k, d_{k-1}, \dots, d_2$ and $\eps_k \ll \eps' \ll \psi, d_k, \beta, 1/k \leq 1/3$ and $1/n \ll 1/t$ such that~$t$ divides~$n$ and $1/d_i \in \mathbb{N}$ for all $2 \leq i \leq k-1$.
Let~$G$ be a $k$-graph on~$n$ vertices and $\mathcal{J}$ be a $(\cdot, \cdot,\eps, \eps_k, r)$-regular slice for~$G$. Further, let $\mathcal{J}$ have~$t$ clusters $V_1, \dots, V_t$ all of size~$n/t$ and density vector $\mathbf{d} = (d_{k-1}, \dots, d_2)$.
Suppose that the reduced graph $\mathcal{R}_{d_k}^{\mathcal{J}}(G)$ contains a tightly connected fractional matching~$\varphi$ with weight~$\mu$. Assume that all edges with non-zero weight have weight at least~$\beta$.
For each $i \in [t]$, let $W_i \subseteq V_i$ be such that $\s{W_i} \geq ((1-3\eps')\varphi(V_i) + \eps')n/t$.
Then $G\left[\bigcup_{i \in [t]}W_i\right]$ contains a tight cycle of length~$\ell$ for each $\ell \leq (1-\psi)k\mu n/t$ that is divisible by~$k$.
\end{lem}

We first explain the main ideas of the proof. We would like to find a regular slice for $G' = G[\bigcup_{i \in [t]} W_i]$ so that we can then apply \cref{lem:cycle} to~$G'$. The issue is that not all vertex classes in~$G'$ have the same size. To get around this we take a refinement of the original partition and use \cref{prop:split} to find a new regular slice with that ground partition. The reduced graph for this new regular slice will be a blow up of the original reduced graph. We can find a corresponding tightly connected matching in this new reduced graph. Then we simply apply \cref{lem:cycle}.

\begin{proof}[Proof of \cref{lem:main}]
Let~$m= n/t$ and $\m = \lfloor \eps'm/2 \rfloor$.
For each $i \in [t]$, let $\widetilde{V}_i \subseteq V_i$ such that $\m \mid \s{\widetilde{V}_i}$ and $\s{V_i \setminus \widetilde{V}_i} \leq \eps'm/2$.
By \cref{regrestriction}, $\mathcal{J}[\widetilde{V}_1, \dots, \widetilde{V}_t]$ is $(\cdot, \cdot, \sqrt{\eps})$-equitable with density vector $(d_{k-1}, \dots, d_2)$.
Let $N = \lfloor m/ \m\rfloor$ and, for each $i \in [t]$, let $N_i = \lfloor ((1- 3\eps')\varphi(V_i) + \eps')N\rfloor \leq \lfloor \s{W_i}/\m \rfloor$.
For each $i \in [t]$, let $V_{i,1}, \dots, V_{i,N}$ be an equipartition of $\widetilde{V}_i$ such that $V_{i,1}, \dots, V_{i,N_i} \subseteq W_i$. 
Let $\widetilde{W} = \{V_{i,j} \colon i \in [t], j \in [N_i]\}$ and $\widetilde{t} = \s{\widetilde{W}}$. 
By \cref{prop:split}, there exists a $(\cdot, \cdot, \eps^{1/4})$-equitable $(k-1)$-complex $\mathcal{J}^*$ with density vector $(d_{k-1}, \dots, d_2)$ and ground partition $\{V_{i,j} \colon i \in [t], j \in [N]\}$ such that $\mathcal{J}[\widetilde{V}_1, \dots, \widetilde{V}_t] = \mathcal{J}^*[\widetilde{V}_1, \dots, \widetilde{V}_t]$.
Let $\J = \mathcal{J}_{\widetilde{W}}^*$, that is~$\J$ is the $(k-1)$-complex contained in $\mathcal{J}^*$ induced by the vertex classes in $\widetilde{W}$.

Let $\widetilde{G}$ be the subgraph of $G[\bigcup \widetilde{W}]$ obtained by removing all edges contained in $k$-tuples of density less than~$d_k$ and in irregular $k$-tuples. We show that~$\J$ is a regular slice for $\widetilde{G}$. Let~$X$ be a set of~$k$ clusters of~$\J$. If the~$k$ clusters in~$X$ are all contained in distinct clusters of $\mathcal{J}$ that form a regular $k$-tuple of density at least~$d_k$, then let~$Y$ denote the $k$-set of these clusters. Note that $(G \cup \mathcal{J})[Y]$ is $((d, d_{k-1},\dots,d_2), \eps_k, \eps, r)$-regular, for some $d \geq d_k -\eps_k$, and thus, by \cref{regrestriction}, $(\widetilde{G} \cup \J)[X]$ is $((d,d_{k-1}, \dots, d_2), \sqrt{\eps_k}, \sqrt{\eps}, r)$-regular. Hence $\widetilde{G}$ is $(d,\sqrt{\eps_k},r)$-regular with respect to $(\J_{X^<})^{(k-1)}$.
Note that, for all other $k$-sets of clusters~$X$, the $k$-partite subgraph of $\widetilde{G}$ induced by the clusters in~$X$ is empty. For these $k$-sets of clusters, $\widetilde{G}$ is $(0,\sqrt{\eps_k},r)$-regular with respect to the polyad $(\J_{X^<})^{(k-1)}$.
Thus~$\J$ is a $(\cdot, \cdot, \sqrt{\eps_k}, \eps^{1/4},r)$-regular slice for $\widetilde{G}$.

Note that $\R = \mathcal{R}_{d_k - 2\sqrt{\eps_k}}^{\J}(\widetilde{G})$ is a blow-up of $\mathcal{R}_{d_k}^{\mathcal{J}}(G)$.
Consider the tightly connected fractional matching~$\varphi$ on $\mathcal{R}_{d_k}^{\mathcal{J}}(G)$ with weight~$\mu$. We construct a tightly connected matching on~$\R$ as follows. 
For each $e \in \mathcal{R}_{d_k}^{\mathcal{J}}(G)$, we will pick a matching~$M_e$ in~$\R$ of size $\widetilde{\varphi}(e) = \lfloor (1-3\eps')\varphi(e)N \rfloor$. 
Note that, for each $i \in [t]$, 
\begin{align}
\label{eq:matching}
\sum_{e \ni V_i} \widetilde{\varphi}(e) \leq \lfloor ((1-3\eps')\varphi(V_i) + \eps')N\rfloor = N_i.
\end{align}
For each vertex~$V_i$ in $\mathcal{R}_{d_k}^{\mathcal{J}}(G)$ and each edge $e \in \mathcal{R}_{d_k}^{\mathcal{J}}(G)$ that contains~$V_i$, we choose disjoint sets $I_{i,e} \subseteq [N_i]$ such that $\s{I_{i,e}} = \widetilde{\varphi}(e)$. This is possible by \cref{eq:matching}.
Recall that~$\R$ is a blow-up of $\mathcal{R}_{d_k}^{\mathcal{J}}(G)$. For each edge $e =\{V_{i_1}, V_{i_2}, \dots, V_{i_k}\} \in \mathcal{R}_{d_k}^{\mathcal{J}}(G)$, the subgraph~$\R_e$ of~$\R$ induced by the set of edges $\{\{V_{i_1,j_1}, \dots, V_{i_k,j_k}\} \colon j_1 \in I_{i_1,e}, \dots, j_k \in I_{i_k,e}\}$ is a balanced complete $k$-partite $k$-graph.
Pick a perfect matching~$M_e$ in~$\R_e$. 
Let $M= \bigcup_{e \in \mathcal{R}_{d_k}^{\mathcal{J}}(G)} M_e$.
Note that~$M$ is a matching of size
\begin{align*}
\sum_{e \in \mathcal{R}_{d_k}^{\mathcal{J}}(G)}\widetilde{\varphi}(e) &= \sum_{e \in \mathcal{R}_{d_k}^{\mathcal{J}}(G)} \lfloor (1- 3\eps')\varphi(e)N\rfloor 
\geq \sum_{\substack{e \in \mathcal{R}_{d_k}^{\mathcal{J}}(G) \\ \varphi(e) > 0}}((1- 3\eps')\varphi(e)N-1)
\\ &\geq (1-3\eps')\mu N- \mu/\beta 
= \br{1-3\eps'- \frac{1}{N\beta}}\mu N 
\\ &\geq (1-3\eps'-\eps'/\beta)\mu N 
\geq (1- \sqrt{\eps'})\mu N
\geq (1-2\sqrt{\eps'})\mu \frac{m}{\m}.
\end{align*}
In the second inequality above we used the fact that since~$\varphi$ is a fractional matching with weight~$\mu$ and all edges have weight at least~$\beta$, there are at most $\mu / \beta$ edges of positive weight.
Since~$\R$ is a blow-up of $\mathcal{R}_{d_k}^{\mathcal{J}}(G)$,~$M$ is tightly connected.
We conclude by applying \cref{lem:cycle} with $k,r,n,\widetilde{t},\psi^2, \eps^{1/4}, \sqrt{\eps_k}, d_k-2\sqrt{\eps_k},d_{k-1}, \dots, d_2, \J,\widetilde{G}, \ell$ playing the roles of $k,r,n_0,t,\psi,\eps,\eps_k,d_k,\dots,d_2, \mathcal{J}, G, \ell$.
\end{proof}

For the next result, we need the following definition.

\begin{defn}
Let $\mu_k^s(\beta,\eps, n)$ be the largest~$\mu$ such that every $2$-edge-coloured $(1-\eps, \eps)$-dense $k$-graph on~$n$ vertices contains a fractional matching with weight~$\mu$ such that all edges with non-zero weight have weight at least~$\beta$ and lie in~$s$ monochromatic tight components.
Let $\mu_k^s(\beta) = \liminf_{\eps \to 0} \liminf_{n \to \infty} \mu_k^s(\beta,\eps,n)/n$.
Similarly, let $\mu_k^*(\beta,\eps,n)$ be the largest~$\mu$ such that every $2$-edge-coloured $(1 - \eps, \eps)$-dense $k$-graph on~$n$ vertices contains a fractional matching with weight~$\mu$ such that all edges with non-zero weight have weight at least~$\beta$ and lie in one red and one blue tight component.
Let $\mu_k^*(\beta) = \liminf_{\eps \to 0} \liminf_{n \to \infty} \mu_k^*(\beta,\eps,n)/n$. 
\end{defn}

The following is the crucial result that reduces finding cycles in the original graph to finding tightly connected matchings in the reduced graph.

\begin{cor}
\label{cor:matchings_to_cycles}
Let $1/n \ll \eta, \beta, 1/k, 1/s$ with $k \ge 3$. 
Let~$K$ be a $2$-edge-coloured complete $k$-graph on~$n$ vertices.
Then the following hold.
\begin{enumerate}[label= \upshape{(\roman*)}]
    \item $K$ contains~$s$ vertex-disjoint monochromatic tight cycles covering at least $(\mu_k^s(\beta) - \eta)k n$ vertices,
    \item $K$ contains two vertex-disjoint monochromatic tight cycles of distinct colours covering at least $(\mu_k^*(\beta) - \eta)k n$ vertices, and
    \item $K$ contains a monochromatic tight cycle of length~$\ell$ for any $\ell \leq (\mu_k^1(\beta) - \eta)k n$ divisible by~$k$.
\end{enumerate}
\end{cor}
\begin{proof}
We prove the first statement. The other two statements can be proved similarly (where for the third statement we additionally make use of the fact that \cref{lem:main} also allows us to control the length of the resulting cycle).
Without loss of generality assume that $\eta \leq 1/3$.
Let $d_k = 1/2$ and $1/t_0 \ll \eps_k \ll \eps' \ll \eps \ll \eta, \beta, 1/k, 1/s$. Note that $\mu_k^s(\beta, \eps, t) \ge (\mu_k^s(\beta) - \eta^2)t$ for all $t \ge t_0$. 
We choose functions~$\widetilde{\eps}(\cdot)$ and~$r(\cdot)$ where~$\widetilde{\eps}(\cdot)$ approaches zero sufficiently quickly and~$r(\cdot)$ increases sufficiently quickly such that for any integer $t^* \geq t_0$ and $d_2, \dots, d_{k-1} \geq 1/t^*$ we may apply \cref{lem:main} with~$\widetilde{\eps}(t^*)$ and~$r(t^*)$ playing the roles of~$\eps$ and~$r$, respectively. 
We apply \cref{lem:regular slice} to obtain~$n_0$ and~$t_1$. Let $\widetilde{\eps} = \widetilde{\eps}(t_1)$ and $r=r(t_1)$.
Let $n_1 \geq n_0$ be large enough such that for all $n \geq n_1$ and $d_2, \dots, d_{k-1} \geq 1/t_1$ we may apply \cref{lem:main}.
Let $n_2 = n_1 + t_1!$.
We show that the theorem holds for all $n \geq n_2$.
Let~$K$ be a $2$-edge-coloured complete $k$-graph on~$n$ vertices.
Let $\widetilde{n} \leq n$ be the largest integer such that~$t_1!$ divides~$\widetilde{n}$.
Let~$\widetilde{K}$ be a complete subgraph of~$K$ on~$\widetilde{n}$ vertices. 
Note that $\widetilde{n} \geq n_1$.
By \cref{lem:regular slice}, there exists a $(t_0, t_1, \widetilde{\eps}, \eps_k, r)$-regular slice~$\mathcal{J}$ for both $\widetilde{K}^{\red}$ and $\widetilde{K}^{\blue}$.
Let~$t$ be the number of clusters of $\mathcal{J}$ and let $(d_{k-1}, \dots ,d_2)$ be the density vector of $\mathcal{J}$. Let $\widetilde{H} = \mathcal{R}_{d_k}^\mathcal{J}(\widetilde{K}^{\red}) \cup \mathcal{R}_{d_k}^\mathcal{J}(\widetilde{K}^{\blue})$ be a $2$-edge-coloured $k$-graph such that $\mathcal{R}_{d_k}^\mathcal{J}(\widetilde{K}^{\red}) \setminus \mathcal{R}_{d_k}^\mathcal{J}(\widetilde{K}^{\blue}) \subseteq \widetilde{H}^{\red}$ and $\mathcal{R}_{d_k}^\mathcal{J}(\widetilde{K}^{\blue}) \setminus \mathcal{R}_{d_k}^\mathcal{J}(\widetilde{K}^{\red}) \subseteq \widetilde{H}^{\blue}$. 
By \cref{lem:reduced graph edge count}, we have $\s{\widetilde{H}} \geq (1-2\eps_k)\binom{t}{k}$. 
By \cref{prop:dense}, there exists a $(1-(2\eps_k)^{1/(4k^2+1)}, (2\eps_k)^{1/(4k^2+1)})$-dense subgraph $H \subseteq \widetilde{H}$ with $V(H) = V(\widetilde{H})$. Since $\eps_k \ll \eps$,~$H$ is $(1- \eps, \eps)$-dense. 
Let~$\varphi$ be a fractional matching in~$H$ of weight $\mu = \mu_k^s(\beta,\eps,t) \geq (\mu_k^s(\beta) - 2 \eta^2)t$ such that all edges with non-zero weight have weight at least~$\beta$ and lie in~$s$ monochromatic tight components $K_1, \dots, K_s$ of~$H$. For each $j \in [s]$, we define a fractional matching~$\varphi_j$ in~$H$ by setting $\varphi_j(e) = \varphi(e)$ if $e \in K_i$ and $\varphi(e) = 0$ otherwise. For each $j \in [s]$, let~$\mu_j$ be the weight of~$\varphi_j$. It follows that $\sum_{j \in [s]} \mu_j = \mu$.

Let $V_1, \dots, V_t$ be the clusters of $\mathcal{J}$.
For each $i \in [t]$ and $j \in [s]$, we define 
\[
w_{i,j} = \max\{\sum_{\substack{e \in H \\ V_i \in e}} \varphi_j(e) - s \eps', \eps'\}.
\]
For each $i \in [t]$, let~$V_{i,1}, \dots, V_{i,s}$ be disjoint subsets of~$V_i$ such that $\s{V_{i,j}} = \lceil w_{i,j}n/t \rceil$. 
By \cref{lem:main}, there exist tight cycles $C_1, \dots, C_s$ in~$K$ such that, for all $j \in [s]$, $\s{C_j} = (1-\eta^2) \mu_j k \widetilde{n}/t$, $C_j \subseteq K\left[\bigcup_{i \in [t]}V_{i,j}\right]$ and~$C_j$ has the same colour as~$K_j$.
Hence~$C_1, \dots, C_s$ are vertex-disjoint and together cover 
\begin{align*}
(1- \eta^2) \mu k \widetilde{n}/t \geq (1-\eta^2)(\mu_k^s(\beta) - \eta^2) k \widetilde{n}
\geq (\mu_k^s(\beta) - \eta) k n
\end{align*}
vertices of~$K$.
\end{proof}

\section{Blueprints}
\label{section:blueprints}
Let~$H$ be a $2$-edge-coloured $k$-graph. 
We define what we call a \emph{blueprint for~$H$} which is an auxiliary graph that can be used as a guide when finding connected matchings in~$H$.
A form of the notion of a blueprint for $k =3$ already appeared in~\cite{Haxell2009}.

\begin{defn}
\label{defn:blueprint}
Let~$\eps >0$, $k \geq 3$ and let~$H$ be a $2$-edge-coloured $k$-graph on~$n$ vertices. We say that a $2$-edge-coloured $(k-2)$-graph~$G$ with $V(G) \subseteq V(H)$ is an \emph{$\eps$-blueprint for~$H$}, if
\begin{enumerate}[label = {$(\text{BP}\arabic*)$}, leftmargin= \widthof{BP1000}]
    \item \label{BP1} for every edge~$e \in G$, there exists a monochromatic tight component~$H(e)$ in~$H$ such that~$H(e)$ has the same colour as~$e$ and $d_{\partial H(e)}(e) \geq (1-\eps)n$ and
    \item \label{BP2} for $e, e' \in G$ of the same colour with $\s{e \cap e'} = k-3$, we have $H(e) = H(e')$.
\end{enumerate}
We say that \emph{$e$ induces~$H(e)$} and write~$R(e)$ or~$B(e)$ instead of~$H(e)$ if~$e$ is red or blue, respectively.
We simply say that~$G$ is a \emph{blueprint}, when~$H$ and~$\eps$ are clear from context. For $S \in \binom{V(H)}{k-3}$, all the red (blue) edges of a blueprint containing~$S$ induce the same red (blue) tight component, so we call that component the red (blue) tight component induced by~$S$. Note that any subgraph of a blueprint is also a blueprint.
\end{defn}

\begin{ex}
\label{ex:blueprint}
Let $k \geq 3$ and let~$n$ be a positive integer. Let~$A$ and~$B$ be disjoint vertex sets with $\s{A \cup B} = n$. Let~$K^{(k)}(A,B)$ be the $2$-edge-coloured complete $k$-graph with vertex set $A \cup B$ where an edge~$e$ is red if and only if $\s{e \cap A}$ is even (and blue otherwise). Let~$H$ be~$K^{(k)}(A,B)$ and let~$G$ be $K^{(k-2)}(A,B)$ with colours reversed. If $\eps \geq \frac{k-2}{n}$, then~$G$ is an $\eps$-blueprint for~$H$. Indeed, for an edge $e \in G$ we can set $H(e) = \{f \in H \colon \s{f \cap A} = \s{e \cap A} +1\}$.
\end{ex}

The main aim of this section is to prove the following lemma that establishes the existence of blueprints for $2$-edge-coloured $(1-\eps, \alpha)$-dense graphs.
\begin{lem}
\label{lem:generalblueprint}
Let $1/n \ll \eps \leq \alpha \ll 1/k \leq 1/3$. Let~$H$ be a $2$-edge-coloured $(1-\eps, \alpha)$-dense $k$-graph on~$n$ vertices.
Then there exists a $3\sqrt{\eps}$-blueprint~$G_*$ for~$H$ with $V(G_*) = V(H)$ and $\s{G_*} \geq (1-\alpha -24k\sqrt{\eps})\binom{n}{k-2}$.
Moreover, if $k \geq 4$ and $\eps \ll \alpha$, there exists a $(1-\alpha^{1/(4(k-2)^2+1)}, \alpha^{1/(4(k-2)^2+1)})$-dense spanning subgraph~$G$ of~$G_*$.
\end{lem}

We need a few simple preliminary results to prove \cref{lem:generalblueprint}. For a $2$-graph $G$, we denote by $\delta(G)$ the minimum degree of $G$.
First we show that any $2$-edge-coloured $2$-graph with large minimum degree contains a large monochromatic connected subgraph. This proposition is implied by \cite[Lemma 1.5]{Gyarfas2012} but we include a proof for completeness. 
\begin{prop}
\label{prop:dcomp}
Let $0 < \beta \leq 1/6$ and let~$F$ be a $2$-edge-coloured $2$-graph with $\s{ V(F) } \leq n$ and $ \delta(F) \geq~(1-\beta)n$. Then there exists a subgraph~$F'$ of~$F$ of order at least~$(1-\beta)n$ that contains a spanning monochromatic component and $\delta(F') \geq (1-2\beta)n$.
\end{prop}
\begin{proof}
Let~$F'$ be an induced subgraph of~$F$ of maximum order that contains a spanning monochromatic component. Assume without loss of generality that~$F'$ contains a spanning red component.
Let~$S = V(F')$ and $\overline{S} = V(F) \setminus V(F')$. 
Since $\delta(F) \geq (1-\beta)n$, we have that $ \s{ S } \geq (1-\beta)n/2 $. 
Suppose, for a contradiction, that $ \s{ S } < (1-\beta)n$. 
Note that all edges between~$S$ and $\overline{S}$ are blue. 
If $\delta(F)-\s{ S } + 1 > \s{ \overline{S}}/2$, then each pair of vertices in~$S$ has a common neighbour in~$\overline{S}$ and so there is a blue component strictly containing~$S$ which contradicts the maximality of~$F'$.
Therefore \[ \delta(F)-\s{ S } + 1 \leq \s{ \overline{S} }/2 = (\s{ V(F) } - \s{ S })/2 \leq (n - \s{ S })/2.\] Hence \[\s{ S } \geq 2 \delta(F) - n + 2 \geq 2(1-\beta)n -n + 2 = (1-2\beta)n+2.\] 
But now every pair of vertices in~$\overline{S}$ has a common neighbour in~$S$, since $ \s{ \overline{S} } = \s{ V(F) } - \s{ S } \leq 2 \beta n $ and so \[\delta(F) - \s{ \overline{S} } + 1 \geq (1-\beta)n - 2\beta n + 1 = (1-3 \beta)n +1 > n/2.\] 
Thus $\overline{S} \cup N_F(\overline{S})$ is spanned by a blue component. But since \[ \s{ \overline{S} \cup N_F(\overline{S}) } \geq \delta(F) \geq (1-\beta)n, \] we have a contradiction.
It is easy to see that $\delta(F') \geq (1-2\beta)n$.
\end{proof}


\begin{prop}
\label{prop:mindeg}
Let $1/n \ll \gamma \leq 1/9$. Let~$F$ be a $2$-graph with $\s{ V(F) } \leq n$ and $\s{ E(F) } \geq (1-\gamma)\binom{n}{2}$. Then there exists a subgraph of~$F$ with minimum degree at least $(1-3\sqrt{\gamma})n$.
\end{prop}
\begin{proof}
Let $W = \{v \in V(F)\colon d(v) < (1- 2\sqrt{\gamma})n\}$. 
We have that \[(1-2\gamma)n^2 \leq 2 \s{E(F)} = \sum_{v \in V(F)}d(v) \leq n^2 - 2\sqrt{\gamma}n\s{W}.\] 
This implies that $\s{ W } \leq \sqrt{\gamma}n $. 
Let~$F^* = F - W$. 
It follows that $\delta(F^*)\geq (1-2\sqrt{\gamma})n-\s{ W} \geq (1-3\sqrt{\gamma})n. $
\end{proof}


\begin{cor}
\label{cor:Ecomp}
Let $1/n \ll \eps \leq 1/324$. Let~$F$ be a $2$-edge-coloured $2$-graph with $\s{ V(F) } \, \leq n$ and $\s{ E(F) } \geq (1-\eps)\binom{n}{2}$. Then there exists a subgraph~$F'$ of~$F$ of order at least $(1-3\sqrt{\eps})n$ that contains a spanning monochromatic component and $\delta(F') \geq (1-6\sqrt{\eps})n$.
\end{cor}
\begin{proof}
By \cref{prop:mindeg}, there exists a subgraph~$F^*$ of~$F$ with $\delta(F^*) \geq (1-3\sqrt{\eps})n$. We conclude by applying \cref{prop:dcomp} with~$F=F^*$ and $\beta = 3\sqrt{\eps}$.
\end{proof}

\subsection{Proof of \texorpdfstring{\cref{lem:generalblueprint}}{the blueprint lemma}}
Now we show that for any $(1-\eps, \alpha)$-dense $2$-edge-coloured graph we can find a dense blueprint.
\begin{proof}[Proof of \cref{lem:generalblueprint}]
Let $F = \partial^2H$. Since~$H$ is $(1-\eps, \alpha)$-dense,
\[
E(F) = \left\{e \in \binom{V(H)}{k-2}\colon d_H(e)>0\right\} = \left\{e \in \binom{V(H)}{k-2}\colon d_H(e)\geq (1-\eps)\binom{n}{2}\right\}
\]
and 
\begin{align}
\label{eq:edges}
\s{E(F)} \geq (1-\alpha)\binom{n}{k-2}.
\end{align}
We now colour each edge~$e$ of~$F$ as follows. Note that the link graph~$H_e$ is a $2$-graph. We induce a 2-edge-colouring on~$H_e$ by colouring the $2$-edge~$f \in H_e$ with the colour of the $k$-edge $e\cup f \in H$. By \cref{cor:Ecomp}, there exists a monochromatic component in~$H_e$ of order at least $(1-3\sqrt{\eps})n$. Let~$K_e$ be such a component chosen arbitrarily. We colour the edge~$e$ according to the colour of~$K_e$. 
If~$e$ is red in~$F$, then we define $R(e) \subseteq H$ to be the red tight component containing all the edges~$e \cup f$ where $f \in K_e$. If~$e$ is blue in~$F$, then we define~$B(e)$ analogously. 

In the next claim we show that, for each $S \in \binom{V(H)}{k-3}$, almost all edges in~$F$ of the same colour containing~$S$ induce the same monochromatic tight component in~$H$.

\begin{claim}
For each $S \in \binom{V(H)}{k-3}$, there exist $\Gamma^{\red}(S) \subseteq N_F^{\red}(S)$ and $\Gamma^{\blue}(S) \subseteq N_F^{\blue}(S)$ with $\s{\Gamma^{\red}(S)} \geq \s{N_F^{\red}(S)} - 6 \sqrt{\eps}n$ and $\s{\Gamma^{\blue}(S)} \geq \s{N_F^{\blue}(S)} - 6 \sqrt{\eps}n$ such that, for all $y_1, y_2 \in \Gamma^{\red}(S)$, $R(S \cup y_1) = R(S \cup y_2)$ and, for all $y_1',y_2' \in \Gamma^{\blue}(S)$, $B(S \cup y_1') = B(S \cup y_2')$.
\end{claim}
\begin{proofclaim}
We only prove the statement for~$N_F^{\red}(S)$ as the proof of the statement for~$N_F^{\blue}(S)$ is analogous. Assume $\s{N_F^{\red}(S)} > 6 \sqrt{\eps}n$ (or else we simply set $\Gamma^{\red}(S) = \varnothing$). 
Let~$D$ be the directed graph with vertex set~$N_F^{\red}(S)$ and edge set 
\[
E(D) = \left\{y_1y_2 \colon y_1 \in V(K_{S \cup y_2})\right\}.
\]
Note that, for $y_1y_2 \in E(D)$, there exists an edge in $R(S \cup y_2)$ containing $S \cup y_1y_2$. So if~$y_1y_2$ is a double edge (that is, $y_1y_2, y_2y_1 \in E(D)$), then $R(S \cup y_1) = R(S \cup y_2)$. 
For $y \in N_F^{\red}(S)$, 
\[
d_D^-(y) \geq \s{N_F^{\red}(S) \cap V(K_{S \cup y})} \geq \s{N_F^{\red}(S)} -3\sqrt{\eps}n,
\]
since $\s{V(K_{S\cup y})} \geq (1-3\sqrt{\eps})n$.
Hence the number of double edges in~$D$ is at least 
\[
\s{N_F^{\red}(S)}\br{\s{N_F^{\red}(S)}-3\sqrt{\eps}n} - \frac{1}{2}\s{N_F^{\red}(S)}^2
= \frac{1}{2}\s{N_F^{\red}(S)}\br{\s{N_F^{\red}(S)}-6\sqrt{\eps}n}.
\]
Thus there exists a vertex $y_0 \in N_F^{\red}(S)$ that is incident to at least $\s{N_F^{\red}(S)} -6\sqrt{\eps}n$ double edges.
Let $\Gamma^{\red}(S) = \{y_0\} \cup \{y \in N_F^{\red}(S)\colon yy_0, y_0y \in E(D)\}$.
Note that $\s{\Gamma^{\red}(S)} \geq \s{N_F^{\red}(S)} - 6\sqrt{\eps}n$ and $R(S \cup y) = R(S \cup y_0)$ for all $y \in \Gamma^{\red}(S)$.
\end{proofclaim}
Consider the multi-$(k-2)$-graph~$D^*$ with 
\[
E(D^*) = \left\{S \cup y \colon S \in \binom{V(H)}{k-3}, y \in \Gamma^{\red}(S)\cup \Gamma^{\blue}(S) \right\}.
\]
Note that 
\begin{align*}
    \s{E(D^*)} &= \sum_{S \in \binom{V(H)}{k-3}} \s{\Gamma^{\red}(S) \cup \Gamma^{\blue}(S)}
    \geq \sum_{S \in \binom{V(H)}{k-3}} (d_F(S) - 12 \sqrt{\eps}n)
    \\ &\geq (k-2)\s{F}-24k\sqrt{\eps} \binom{n}{k-2}. 
\end{align*}
Observe that every edge of $D^*$ is an edge of $F$.
Every edge in~$D^*$ has multiplicity at most~$k-2$. So at least $\s{F} -24k\sqrt{\eps}\binom{n}{k-2}$ edges $e \in \binom{V(H)}{k-2}$ have multiplicity~$k-2$ in~$D^*$. 
Let~$G_*$ be the $(k-2)$-graph on~$V(H)$ such that $e \in G_*$ if and only if~$e$ has multiplicity~$k-2$ in~$D^*$.
So, by (\ref{eq:edges}), $\s{G_*} \geq \s{F} -24k\sqrt{\eps}\binom{n}{k-2} \geq (1- \alpha -24k\sqrt{\eps})\binom{n}{k-2}$.

We now show that~$G_*$ is a $3\sqrt{\eps}$-blueprint for~$H$.
Consider any $e, e' \in G_*^{\red}$ with $\s{e \cap e'} = k-3$.
Let $S = e\cap e'$, $y = e'\setminus S$ and $y' = e \setminus S$. 
Since $e,e' \in G_*^{\red}$, we have $y,y' \in \Gamma^{\red}(S)$ and so $R(e)=R(S\cup y) =R(S\cup y') = R(e')$. Further, for~$e \in G_*^{\red}$, we have $d_{\partial R(e)}(e) \geq \s{V(K_e)} \geq (1-3\sqrt{\eps})n$. Analogous statements hold for edges of~$G_*^{\blue}$.

If $k \geq 4$ and $\eps \ll \alpha$, then $\s{G_*} \geq (1-2\alpha)\binom{n}{k-2}$ and thus by \cref{prop:dense} there exists a subgraph $G \subseteq G_*$ such that~$G$ is $(1-\alpha^{1/(4(k-2)^2+1)}, \alpha^{1/(4(k-2)^2+1)})$-dense and $V(G) = V(G_*) = V(H)$. 
\end{proof}

\subsection{Some lemmas about blueprints}
Let~$H$ be a $k$-graph and~$G$ be a blueprint for~$H$. 
We write~$H(G)$ for $\bigcup_{e \in G} H(e)$.
We write~$G^+$ for the subgraph of~$H(G)$ with edge set 
\[
E(G^+) = \{e \in H(G) \colon f \subseteq e \text{ for some } f \in G\}, 
\]
that is, the subgraph of~$H(G)$ obtained by deleting all edges that do not contain an edge of~$G$.
Note that this also defines $(G')^+$ for any subgraph $G'$ of $G$ as a subgraph of a blueprint for $H$ is also a blueprint for $H$.
Moreover, note that~$G^+$ is a subgraph of~$H$, not of~$G$. 
For a red tight component~$R_*$ and a blue tight component~$B_*$ in~$H$, we denote by~$R_*^{k-2}$ and~$B_*^{k-2}$ the edges of~$G$ that induce~$R_*$ and~$B_*$, respectively. 

We prove some lemmas that we will use several times later on.
Roughly speaking, the following lemma says that if~$S$ is a set of~$k-4$ vertices of~$H$ contained in many edges of~both~$R_*^{k-2}$ and~$B_*^{k-2}$, then~$S$ is contained in an edge of~$R_*$ or~$B_*$.

\begin{lem}
\label{lem:vertexdeg}
Let $1/n \ll \eps \ll \alpha \ll 1.$
Let~$H$ be a $2$-edge-coloured $(1-\eps, \alpha)$-dense $k$-graph on~$n$ vertices and~$G$ a $3\sqrt{\eps}$-blueprint for~$H$.
Let~$R_*$ and~$B_*$ be a red and a blue tight component of~$H$, respectively.
Let $U \subseteq V(G)$ and $S \in \binom{U}{k-4}$ such that 
\[
d_{R_*^{k-2}}(S,U), d_{B_*^{k-2}}(S,U) \geq \eps^{1/4}n^2.
\]
Then there exist $x,x',y,y' \in U$ such that $S \cup xx' \in R_*^{k-2}$, $S \cup yy' \in B_*^{k-2}$, $S \cup xx'y \in \partial R_*$, $S \cup yy'x \in \partial B_*$ and $S \cup xx'yy' \in H$. In particular, $(R_*^{k-2})^+[U] \cup (B_*^{k-2})^+[U] \neq \varnothing$.
\end{lem}
\begin{proof}
Let $X_{R_*} = \{x \in U \colon d_{R_*^{k-2}}(S \cup x, U) \geq \eps^{1/2}n\}$ and $X_{B_*} = \{x \in U \colon d_{B_*^{k-2}}(S \cup x, U) \geq \eps^{1/2} n\}$. 
Note that 
\[
\eps^{1/4} n^2 \leq d_{R_*^{k-2}}(S, U) = \frac{1}{2} \sum_{x \in U} d_{R_*^{k-2}}(S \cup x, U) \leq n \s{X_{R_*}} + \eps^{1/2}n^2.
\]
Thus $\s{X_{R_*}} \geq (\eps^{1/4} - \eps^{1/2}) n \geq \frac{1}{2}\eps^{1/4} n$. Similarly, $\s{X_{B_*}} \geq \frac{1}{2} \eps^{1/4}n$.

For each $x \in X_{R_*}$, let 
\begin{align*}
Y_x &= \{y \in X_{B_*} \colon S \cup yy' \in B_*^{k-2} \text{ and } S \cup xyy' \in \partial B_* \text{ for some } y' \in U \} \\ 
&= \bigcup_{y' \in U} N_{B_*^{k-2}}(S \cup y') \cap N_{\partial B_*}(S \cup xy').
\end{align*}
For each $y \in X_{B_*}$, there exists $y' \in U$ with $S \cup yy' \in B_*^{k-2}$. By \ref{BP1}, $d_{\partial B_*}(S \cup yy', X_{R_*}) \geq \s{X_{R_*}} - 3\sqrt{\eps}n$. Hence each $y \in X_{B_*}$ is contained in at least $\s{X_{R_*}} - 3 \sqrt{\eps}n$ of the sets~$Y_x$. By averaging, there exists an $x \in X_{R_*}$ such that 
\[
\s{Y_x} \geq
\frac{(\s{X_{R_*}}- 3\sqrt{\eps}n)\s{X_{B_*}}}{2 \s{X_{R_*}}} \geq \frac{1}{4} \s{X_{B_*}} \geq \frac{1}{8} \eps^{1/4} n.
\]
Fix such an $x \in X_{R_*}$. For each $y \in Y_x$, choose a vertex $y' \in U$ such that $S \cup yy' \in B_*^{k-2}$ and $S \cup xyy' \in \partial B_*$.
Let $X = N_{R_*^{k-2}}(S \cup x, U)$, so $\s{X} \geq \eps^{1/2}n$, since $x \in X_{R_*}$.
For each $y \in Y_x$, since~$H$ is $(1-\eps, \alpha)$-dense, there are at least $\s{X} - \eps n$ vertices $x' \in X$ such that $S \cup xx'yy' \in H$. Thus, by averaging, there exists a vertex $x' \in X$ and a set $\widetilde{Y}_x \subseteq Y_x$ with 
\[
\s{\widetilde{Y}_x} \geq \frac{(\s{X} - \eps n) \s{Y_x}}{2 \s{X}} \geq \frac{1}{4} \s{Y_x} \geq \frac{1}{32}\eps^{1/4} n
\]
such that $S \cup xx'yy' \in H$ for all $y \in \widetilde{Y}_x$.
Fix such an $x' \in X$.
Since $S \cup xx' \in R_*^{k-2}$, we have that 
\[
\s{N_{\partial R_*}(S \cup xx') \cap \widetilde{Y}_x} \geq \s{\widetilde{Y}_x} - 3\sqrt{\eps}n \geq \br{\frac{1}{32}\eps^{1/4} -3\sqrt{\eps}}n >0.
\]
Choose $y \in N_{\partial R_*} (S \cup xx') \cap \widetilde{Y}_x$.
We have $S \cup xx' \in R_*^{k-2}$, $S \cup yy' \in B_*^{k-2}$, $S \cup xx'y \in \partial R_*$, $S \cup xyy' \in \partial B_*$ and $S \cup xx'yy' \in H$ as required.
\end{proof}

The following lemma shows that if we have a vertex set~$T \in \binom{V(G)}{k-3}$ such that $d_G^{\red}(T)$ and $d_G^{\blue}(T)$ are both large, then~$T$ is contained a lot of sets in $\partial R \cap \partial B$, where~$R$ and~$B$ are the red and blue tight components induced by the red and blue edges incident to~$T$, respectively.

\begin{lem}
\label{lem:shadow}
Let $1/n \ll \eps \ll 1$, $k \ge 3$ and $\delta > 5 \sqrt{\eps}$.  
Let~$H$ be a $2$-edge-coloured $k$-graph on~$n$ vertices and~$G$ a $3\sqrt{\eps}$-blueprint for~$H$.
Let $T \in \binom{V(H)}{k-3}$.
Let $S^{\blue} \subseteq N_G^{\blue}(T)$ and $S^{\red} \subseteq N_G^{\red}(T)$ be such that $|S^{\blue}|, |S^{\red}| \ge \delta n $. 
Then there exists a vertex $y \in S^{\blue}$ such that, for
\begin{align*}
\Gamma^{\red}_y = \{x \in S^{\red} \colon T \cup xy \in \partial R(T \cup x) \cap \partial B(T \cup y)\},
\end{align*}
we have $\s{\Gamma^{\red}_y} \ge ( \delta - 6 \sqrt{\eps})n $. 
Moreover, if $\delta \geq \eps^{1/9}$, then $\s{\Gamma^{\red}_y}\geq (1-\eps^{1/4}) \s{S^{\red}} $.
The same statements hold when the colours are reversed.
\end{lem}

\begin{proof}
Let $m_{\blue} = |S^{\blue}|$ and $m_{\red} = |S^{\red}|$.
If $\delta < \eps^{1/9}$, then we may assume that $m_{\blue} = m_{\red} = \lceil \delta n \rceil$ by deleting vertices in~$S^{\blue}$ and~$S^{\red}$ if necessary. 
Let~$D$ be the bipartite directed graph with vertex classes~$S^{\blue}$ and~$S^{\red}$ such that, for each $y \in S^{\blue}$ and $x \in S^{\red}$, we have 
$N_D^+(y) = N_{\partial B} (T \cup y) \cap S^{\red}$ and $N_D^+(x) = N_{\partial R}(T \cup x) \cap S^{\blue}$.
Since~$G$ is a $3\sqrt{\eps}$-blueprint for~$H$, we have that 
\begin{align*}
\s{ E(D) } \geq m_{\blue}(m_{\red} - 3\sqrt{\eps}n) + m_{\red}(m_{\blue} - 3\sqrt{\eps}n) = 2m_{\blue}m_{\red} - 3\sqrt{\eps}n(m_{\blue} + m_{\red}).
\end{align*}
Thus the number of double edges in~$D$ is at least $m_{\blue}m_{\red} -3\sqrt{\eps}n(m_{\blue}+m_{\red})$. 
For each $y \in S^{\blue}$, let $\Gamma_y = \{x \in S^{\red} \colon xy, yx \in D\}$.
Hence there is some vertex $y \in S^{\blue}$ such that 
\begin{align*}
    \s{\Gamma_y} \geq m_{\red}-3\sqrt{\eps}n\br{\frac{m_{\blue} + m_{\red}}{m_{\blue}}} \geq
    \begin{cases}
    (\delta - 6 \sqrt{\eps})n, &\text{if } \delta < \eps^{1/9}; 
    \\ m_{\red}(1-\eps^{1/4}), 
    &\text{otherwise.}
    \end{cases}
\end{align*}
Note that if $xy, yx \in D$ with $x \in S^{\red}$ and $y \in S^{\blue}$, then 
$T \cup xy \in \partial R(T \cup x) \cap \partial B(T \cup y)$.
Hence $\Gamma_y \subseteq \Gamma_y^{\red}$ and thus the lemma follows.
\end{proof}

Roughly speaking, in the next lemma we consider the following situation. 
Let~$R$ be a red tight component in~$H$,~$G$ be a blueprint for~$H$ and $R_G \subseteq G^{\red}$ be such that $H(R_G) \subseteq R$.
We pick a maximal matching in~$R_G^+$ and let~$U$ be the remaining vertices of~$H$ not in this matching, so~$R_G^+[U]$ is empty.
Then the lemma implies that the number of monochromatic tight components in~$U$ is less than what we would expect. 
In particular, if $k = 4$, then the edges in~$G[U]$ induce only two monochromatic tight components in~$H$. 

\begin{lem}
\label{lem:reducing_components}
Let $k \ge 4$ and $1/n \ll \eps \ll \alpha, \delta \ll \eta \ll 1$. 
Let~$H$ be a $(1-\eps, \alpha)$-dense $k$-graph and~$G$ a $3\sqrt{\eps}$-blueprint for~$H$. 
Let~$R$ be a red tight component in~$H$.
Let~$R_G \subseteq G^{\red}$ be such that $H(R_G) \subseteq R$. 
Let $U \subseteq V(H)$ be such that $\s{U } \geq \eta n/2$ and $R_G^+[U] = \varnothing$. 
Let $S \in \binom{U}{k-4}$ be such that the link graph~$G_S$ of~$G$ satisfies $G_S^{\red}[U ] \subseteq (R_G)_S$ and $\delta(G_S[U ]) \geq \s{U } - \delta n$. 
Then there exists a subgraph~$J_S$ of~$G_S[U]$ such that $\s{J_S} \geq \s{G_S[U]} - 7\delta^{1/4}n^2$ and $H(S \cup e) = H(S \cup e')$ for all $e,e' \in J_S$ of the same colour. 
In particular, if $k=4$, then the edges in~$J$ induce only one red and one blue tight component in~$H$.
The same statement holds when the colours are reversed.
\end{lem}

\begin{proof}
Set $J_S^{\red} = G_S^{\red}[U ]$.
Note that for $e, e' \in J_S^{\red}$, we have $e, e' \in (R_G)_S$ and thus $H(S \cup e) = H(S \cup e') = R$ since $H(R_G) \subseteq R$.
Therefore to prove the lemma, it suffices to prove that there exists $J_S^{\blue} \subseteq G_S^{\blue}[U ]$ such that $\s{J_S^{\red}}+\s{J_S^{\blue}} \geq \s{G_S[U ]} - 7\delta^{1/4}n^2$ and $H(S \cup e) = H(S \cup e')$ for all $e,e' \in J_S^{\blue}$.

For simplicity we assume $k=4$ and $S = \varnothing$. It is easy to see that an analogous argument works in the general case. Thus for the rest of the proof, we omit the subscript~$S$.

Let $K = G[U ]$.
If $\s{K^{\blue}} < 2 \delta^{1/2}n^2$, then we are done by setting $J^{\blue} = \varnothing$ as 
\begin{align*}
    \s{J^{\red}} = \s{K^{\red}} = \s{K} - \s{K^{\blue}} \geq \s{K} - 2\delta^{1/2}n^2 \geq \s{K} - 7\delta^{1/4}n^2.
\end{align*}
Now assume $\s{K^{\blue}} \geq 2 \delta^{1/2}n^2$.
Let $X = \{x \in V(K) \colon d_K^{\blue}(x) \geq \delta n\}$.
We have that 
\begin{align*}
    2\delta^{1/2} n^2 \leq \s{K^{\blue}} &\leq \sum_{x \in U } d_K^{\blue}(x) \leq n\s{X} + \delta n^2.
\end{align*}
Thus $\s{X} \geq \delta^{1/2} n$.
Let~${D}$ be the digraph with vertex set~$X$ such that, for each $x \in X$,
\begin{align*}
	N_{D}^+(x) & = N_{K}^{\blue}(x, X) \cup \{ x' \in N_{K}^{\red}(x,X) \colon xx'y \in \partial R \cap \partial B(xy) \text{ for some } y \in N_{K}^{\blue}(x)\}.
\end{align*}
We now bound $\delta^+(D)$ as follows. 
If $d_K^{\red}(x,X) \geq \delta n$, then by applying \cref{lem:shadow} (with $x, N_G^{\blue}(x, U ), N_G^{\red}(x,X), \delta$ playing the roles of $T,S^{\blue},S^{\red}, \delta$), we deduce that 
\begin{align*}
	&\s{\{x' \in N_{K}^{\red}(x,X) \colon xx'y \in \partial R(xx') \cap \partial B(xy) \text{ for some } y \in N_{K}^{\blue}(x)\}} \\ &\ge (1-\eps^{1/4})d_K^{\red}(x,X).
\end{align*}
Recall that $R = R(xx')$ for all $x' \in N_{K}^{\red}(x,X)$, $\s{X} \geq \delta^{1/2}n$ and $\eps \ll \delta$.
Hence 
\begin{align*}
    d_{D}^+(x) &\geq d_K^{\blue}(x,X) + (1-\eps^{1/4})d_K^{\red}(x,X) 
    \geq (1-\eps^{1/4})(d_K^{\blue}(x,X) + d_K^{\red}(x,X)) \\
		& = (1-\eps^{1/4})d_K(x,X)
    \geq (1-\eps^{1/4})(\s{X} - \delta n) 
    \geq (1-2\delta^{1/2})\s{X}.
\end{align*}
On the other hand, if $d_K^{\red}(x,X) < \delta n$, then 
\begin{align*}
    d_{D}^+(x)\geq d_K^{\blue}(x,X) \geq \s{X} - \delta n - d_K^{\red}(x,X)
    \geq \s{X} - 2\delta n
		\geq (1-2 \delta^{1/2})\s{X}.
\end{align*}
Therefore, we have $\delta^+(D) \geq (1-2\delta^{1/2})\s{X}$ and so 
$\s{E({D})} \geq (1-2\delta^{1/2}) \s{X}^2
\geq 2(1-2\delta^{1/2})\binom{\s{X}}{2}$.
Let~$F $ be the graph with vertex set~$X$ in which~$xx'$ forms an edge if and only if it forms a double edge in~${D}$.
Note that $|F| \ge (1-4\delta^{1/2}) \binom{\s{X}}{2}$.
By \cref{prop:mindeg}, there exists a subgraph~$F ^*$ of~$F $ with $\delta(F ^*) \geq (1 -6\delta^{1/4})\s{X}$.
Clearly, $F ^*$ is connected.

Let $J^{\blue} = \{ xx' \in K^{\blue} \colon x \in V(F ^*)\}$.
We have
\begin{align*}
    \s{J^{\red} \cup J^{\blue}} &\geq \s{K} - \sum_{x' \in U \setminus X}d^{\blue}_K(x') - \s{X\setminus V(F ^*)} n\\
		& \geq 		\s{K} - \delta n^2 - 6 \delta^{1/4} n^2 
		\geq \s{G[U ]} - 7\delta^{1/4}n^2.
\end{align*}
We now show that $B(x_1z_1 ) = B(x_2z_2 )$ for all $x_1z_1, x_2z_2 \in J^{\blue}$.
Since~$F ^*$ is connected and $d_{J^{\blue}}(x) > 0$ for all $x \in V(F^*)$, it suffices to consider the case when $x_1 x_2 \in F ^*$.
If $x_1x_2 \in K^{\blue}$, then $x_1z_1 , x_1x_2 , x_2z_2   \in G^{\blue}$ and so $B(x_1z_1 ) = B(x_1x_2 ) = B(x_2z_2 )$, since~$G$ is a blueprint.
Now assume that $x_1x_2 \in K^{\red}$. 
Since $x_1x_2 \in F ^* \subseteq F $, there are $y_1 \in N_K^{\blue}(x_1)$ and $ y_2 \in N_K^{\blue}(x_2)$ such that $x_1x_2y_1   \in \partial R \cap \partial B(x_1y_1 )$ and $x_1x_2y_2   \in \partial R \cap \partial B(x_2y_2 )$.
Let $u \in N_H(x_1x_2y_1 ) \cap N_H(x_1x_2y_2 ) \cap U$.
Since $R^+_G[U] = \varnothing$, we have $x_1x_2y_1u , x_1x_2y_2u   \in H^{\blue}$. 
Hence, $B(x_1y_1 ) = B(x_2y_2 )$.
Moreover, since $x_1y_1 , x_1z_1 , x_2y_2 , x_2z_2   \in G^{\blue}$, we have $B(x_1z_1 ) = B(x_1y_1 )= B(x_2y_2 ) = B(x_2z_2 )$ as required. 
\end{proof}

\section{Monochromatic connected matchings in \texorpdfstring{$K_n^{(4)}$}{Kn(4)}}
\label{section:matchings Kn4}

In this section, we prove that every almost complete red-blue edge-coloured $4$-graph~$H$ contains a red and a blue tightly connected matching that are vertex-disjoint and together cover almost all vertices of~$H$.

\begin{lem}
\label{lem:matchings}
Let $1/n \ll \eps \ll \alpha \ll \eta < 1$.
Let~$H$ be a $2$-edge-coloured $(1-\eps,\alpha)$-dense $4$-graph on~$n$ vertices. Then~$H$ contains two vertex-disjoint monochromatic tightly connected matchings of distinct colours such that their union covers all but at most~$3\eta n$ of the vertices of~$H$.
\end{lem}

Note that this implies $\mu_4^*(1, \eps, n) \geq (1- 3\eta)n/4$ for $1/n \ll \eps \ll \eta < 1$. Hence $\mu_4^*(1) \geq 1/4$.
Therefore, together with \cref{cor:matchings_to_cycles}, \cref{lem:matchings} implies \cref{thm:1}.

To prove \cref{lem:matchings} we first need the following lemma which chooses the initial tight components in~$H$ in which we find our tightly connected matchings. 

\begin{lem}
\label{lem:4existence}
Let $1/n \ll \eps \ll \alpha \ll \eta < 1$.
Let~$H$ be a $2$-edge-coloured $(1-\eps,\alpha)$-dense $4$-graph on~$n$ vertices. 
Suppose that~$H$ does not contain two vertex-disjoint monochromatic tightly connected matchings of distinct colours such that their union covers all but at most~$3\eta n$ of the vertices of~$H$.
Then, there exists a red tight component~$R$ in~$H$, a blue tight component~$B$ in~$H$, 
a $3\sqrt{\eps}$-blueprint~$G$ for~$H$ with $\delta(G) \geq (1-\alpha^{1/30})n$ and a matching~$M_0$ in $R \cup B$ such that the following holds, where $W_0 = V(G)\setminus V(M_0) $.
\begin{enumerate}[label = \upshape (\roman*)]
    \item $ R(e) =R$ and $ B(e') =B$ for all edges $e \in G^{\red}[V(M_0^{\red})\cup W_0]$ and all edges $e' \in G^{\blue}[V(M_0^{\blue})\cup W_0]$, \label{itm:4-3}
    \item $M_0 \subseteq (G^{\red})^+ \cup (G^{\blue})^+$, \label{itm:4-4}
	\item $(G^{\red})^+[W_0] \cup (G^{\blue})^+[W_0]$ is empty. \label{itm:4-5}
\end{enumerate}
\end{lem}
\begin{proof}
By \cref{lem:generalblueprint}, there exists a $3\sqrt{\eps}$-blueprint~$G_0$ for~$H$ with $V(G_0) = V(H)$ and $\s{G_0} \geq (1-\alpha - 96\sqrt{\eps}) \binom{n}{2} \geq (1-4\alpha)\binom{n}{2}$. 
By \cref{cor:Ecomp}, there exists a subgraph~$G_1$ of~$G_0$ of order at least $(1-6\sqrt{\alpha})n$ that contains a spanning monochromatic component and $\delta(G_1) \geq (1-12\sqrt{\alpha})n$. 
Note that that~$G_1$ is also a $3\sqrt{\eps}$-blueprint for~$H$.

We assume without loss of generality that~$G_1$ contains a spanning red component. 
Since~$G_1$ is a blueprint, all the red edges in~$G_1$ induce the same red tight component~$R$ in~$H$.
Let $R^+ = (G^{\red}_1)^+ \subseteq R$. 
Let~$M$ be a matching in~$R^+$ of maximum size. 
Let $U = V(G_1) \setminus V(M)$.

Thus $\s{ U } \geq \eta n$ (or else $\s{V(M)} \geq |V(G_1)| - \s{U} \ge (1-2\eta)n$, a contradiction).
Moreover, $R^+[U] = \varnothing$. 
Since $\delta(G_1) \geq (1-12\sqrt{\alpha})n$, we have $\delta(G_1[U]) \geq \s{U} -\alpha^{1/3}n$.
Hence, by \cref{lem:reducing_components} (with $4,U, \varnothing, \alpha^{1/3}$ playing the roles of $k,U,S,\delta$), there exists a subgraph~$J$ of~$G_1[U]$ such that $\s{J} \geq \s{G_1[U]} - 2\alpha^{1/13}n^2$, such that $H(e) = H(e')$ for all $e, e' \in J$ of the same colour. 
Let $G_2 = (G_1- G_1^{\blue}[U]) \cup J$ and $B = B(e)$ for $e \in J^{\blue}$.
Note that $\s{G_2} \geq (1-\alpha^{1/14})\binom{n}{2}$. 
By \cref{prop:mindeg}, there exists a subgraph~$G$ of~$G_2$ such that $\delta(G) \geq (1 - \alpha^{1/30})n$. 

Let $W = V(G) \setminus V(M)$.
Next, we show that \ref{itm:4-3} and \ref{itm:4-4} hold but with~$M,W$ instead of~$M_0,W_0$.
Note that $M^{\blue} = \varnothing$, so \ref{itm:4-4} holds by our construction.
Since $G^{\red} \subseteq G_1^{\red}$ and~$G_1^{\red}$ is connected and a blueprint, $R(e) = R$ for all $e \in G^{\red}$.
Note that $G^{\blue}[V(M^{\blue})\cup W] = G^{\blue} - V(M) \subseteq G_2^{\blue}[U] = J^{\blue}$, so $B(e) = B$ for all $e \in G^{\blue}[V(M^{\blue})\cup W]$.
Hence \ref{itm:4-3} holds. 
We now add vertex-disjoint edges of $(G^{\red})^+[W] \cup (G^{\blue})^+[W]$ to~$M$ and call the resulting matching~$M_0$. 
We deduce that~$M_0$ satisfies \ref{itm:4-3}--\ref{itm:4-5}.
\end{proof}

We now prove \cref{lem:matchings}. 

\begin{proof}[Proof of \cref{lem:matchings}]
Suppose the contrary that~$H$ does not contain two vertex-disjoint monochromatic tightly connected matchings of distinct colours such that their union covers all but at most~$3\eta n$ of the vertices of~$H$.
We call this the initial assumption.
Apply \cref{lem:4existence} and obtain a red tight component~$R$, a blue tight component~$B$ in~$H$, a $3\sqrt{\eps}$-blueprint~$G$ for~$H$ with $\delta(G) \geq (1-\alpha^{1/30})n$ and a matching~$M_0$ in $R \cup B$ 
satisfying \cref{lem:4existence}\ref{itm:4-3}--\ref{itm:4-5}.

We now fix~$G$,~$R$ and~$B$.
We use the following notation for the rest of the proof. 
For a matching~$M$ in $R \cup B$, we set 
\begin{align*}
W & = W(M) = V(G)\setminus V(M),\\
W_\textup{red} &= W_\textup{red}(M) = \{w \in W \colon d_{G[W]}^{\blue}(w) \leq 8\sqrt{\eps}n\},\\
W_\textup{blue} & = W_\textup{blue}(M) = \{w \in W\colon d_{G[W]}^{\red}(w) \leq 8\sqrt{\eps}n\}.
\end{align*}
Note that $\s{W}\geq \eta n$ by the initial assumption.
Without loss of generality, $\s{W_\textup{blue} (M_0)} \le \s{W_\textup{red} (M_0)}$.

We define $\mathcal{M}$ be the set of matchings~$M$ in $R \cup B$ such that 
\begin{enumerate}[label = \upshape (\roman*$'$)]
		 \item \label{itm:M1} $ R(e) =R$ and $ B(e') =B$ for all edges $e \in G^{\red}[W]$ and $e' \in G^{\blue}[V(M^{\blue})\cup W]$,
		\item $M^{\blue} \subseteq (G^{\blue})^+$, \label{itm:M2}
		\item $(G^{\red})^+[W] \cup (G^{\blue})^+[W]$ is empty. \label{itm:M3}
\end{enumerate}
Note that \ref{itm:M1} and \ref{itm:M2} are weaker statements of those in \cref{lem:4existence}\ref{itm:4-3} and~\ref{itm:4-4}, so $M_0 \in \mathcal{M}$. 
Let $\mathcal{M}'$ be the set of $M \in \mathcal{M}$ also satisfying
\begin{enumerate}[label = \upshape (\roman*$'$),resume]
		\item $\s{W_\textup{blue}} \le \s{W_\textup{red}}$.\label{itm:M4}
\end{enumerate}
Observe that $M_0 \in \mathcal{M}'$, so $\mathcal{M}'$ is nonempty.

Let $\gamma = 10\alpha^{1/30}$.
We now show that, for all $M \in \mathcal{M}$, $W_\textup{red}$ and~$W_\textup{blue}$ partition~$W$, and moreover one of them is small.

\begin{claim}
\label{claim:mathcalM}
Let $M \in \mathcal{M}$. 
The following holds:
\begin{enumerate}[label = \upshape (\alph*)]
		 \item \label{itm:Ma} for all $w \in W$, either $d_{G[W]}^{\red}(w) \le 7\sqrt{\eps}n$ or $d_{G[W]}^{\blue}(w) \le 7\sqrt{\eps}n$,
		\item \label{itm:Mb}$W_\textup{red}$ and~$W_\textup{blue}$ partition~$W$,
		\item \label{itm:Mc} either $\s{ W_\textup{blue}} \le \gamma n $ or $ \s{ W_\textup{red}} \leq \gamma n$.
\end{enumerate}
In particular, if $M \in \mathcal{M}'$, then $\s{W_\textup{blue}} \le \gamma n $.
\end{claim}

\begin{proofclaim}
Suppose that $w \in W$ satisfies $d_{G[W]}^{\red}(w), d_{G[W]}^{\blue}(w) > 7\sqrt{\eps}n$. 
By \cref{lem:shadow} (with $7\sqrt{\eps},w, N_{G[W]}^{\red}(w),N_{G[W]}^{\blue}(w)$ playing the roles of $\delta, T, S^{\red}, S^{\blue}$), there exist $x \in N_{G[W]}^{\red}(w)$ and $y \in N_{G[W]}^{\blue}(w)$ such that $wxy \in \partial R \cap \partial B$. 
In particular, $d_H(wxy) \neq 0$ and thus $d_H(wxy) \geq (1-\eps)n$, which implies that there exists a vertex $w' \in W$ such that $ww'xy \in H$.
Note that $ww'xy \in (G^{\red})^+[W] \cup (G^{\blue})^+[W]$ contradicting~\ref{itm:M3}. 
Hence, $\min \{ d_{G[W]}^{\red}(w), d_{G[W]}^{\blue}(w)\} \le 7\sqrt{\eps}n$.
Since $\delta(G) \geq (1-\alpha^{1/30})n$ implies that $\delta(G[W]) \ge |W| - \alpha^{1/30} n > 16 \sqrt{\eps}n$, we deduce that \ref{itm:Ma} and \ref{itm:Mb} hold.

Recall that $|W| \ge \eta n > 2 \gamma n$. 
So one of $W_\textup{red}$ and $W_\textup{blue}$ has size greater than $\gamma n$. 
Suppose both are (that is, \ref{itm:Mc} is false). 
Since $\delta(G) \geq (1-\alpha^{1/30})n = (1-\gamma/10)n$, we have that there are at least
\begin{align*}
\s{W_\textup{blue}} \br{\s{ W_\textup{red}} -\gamma n/10 - 8\sqrt{\eps}n} \geq \s{W_\textup{blue}} \br{\s{ W_\textup{red} } - \gamma n/5} > 3\s{W_\textup{red}}\s{W_\textup{blue}}/4
\end{align*}
blue edges between~$W_\textup{blue}$ and~$W_\textup{red}$ and similarly there are at least $3\s{ W_\textup{red}} \s{ W_\textup{blue}} /4$ red edges between~$W_\textup{blue}$ and~$W_\textup{red}$.
Thus $e(W_\textup{red},W_\textup{blue}) > \s{W_\textup{red}}\s{W_\textup{blue}}$, a contradiction.
\end{proofclaim}

Let $M_* \in \mathcal{M}'$ be such that $(\s{M_*}, \s{M_*^{\red}})$ is lexicographically maximum. 
We write $W^*, W_\textup{red}^*,W_\textup{blue}^*$ for $W(M_*), W_\textup{red}(M_*),W_\textup{blue}(M_*)$, respectively. 

The next claim shows that almost all $4$-edges in~$H[W^*]$ are blue and they form a tight component. 
Indeed, this follows from the fact that almost all edges in~$G[W^*]$ are red and thus almost all triples in~$W^*$ are in~$\partial R$.

\begin{claim}
\label{claim:B'}
There exists a blue tight component~$B'$ in~$H$ such that the number of triples $xyz \in \binom{W_\textup{red}^*}{3} \cap \partial B'$ with $d_{B'}(xyz, W_\textup{red}^*) \geq \s{W_\textup{red}^*} - \eps n$ is at least $(1-\alpha^{1/31})\s{\binom{W_\textup{red}^*}{3}}$. 
\end{claim}

\begin{proofclaim}
Let~$\mathcal{T}$ be the set of triples $xyz \in \binom{W_\textup{red}^*}{3} \cap \partial R$ such that $xy \in G^{\red}$.
Note that, for any $x \in W_\textup{red}^*$, $y \in N^{\red}_{G}(x, W_\textup{red}^*)$ and $z \in N_{\partial R}(xy, W_\textup{red}^*)$, we have $xyz \in \mathcal{T}$. Thus 
\begin{align*}
    \s{\mathcal{T}} &\geq \frac{1}{3!}\s{W_\textup{red}^*}\br{\s{W_\textup{red}^*}-\alpha^{1/30}n-8\sqrt{\eps}n}\br{\s{W_\textup{red}^*}-3\sqrt{\eps}n} \\
    &\geq \frac{\s{W_\textup{red}^*}^3}{3!}\br{1-\frac{2\alpha^{1/30}n}{\s{W_\textup{red}^*}}} \geq \br{1-\alpha^{1/31}}\s{\binom{W_\textup{red}^*}{3}},
\end{align*}
as $\s{W_\textup{red}^*} \geq \eta n/2$.
By \ref{itm:M3}, we have that if $xyz \in \mathcal{T}$ and $w \in N_H(xyz, W_\textup{red}^*)$, then $wxyz \in H^{\blue}$. 
For $xyz \in \mathcal{T}$, let~$B(xyz)$ be the maximal blue tight component containing all the edges~$xyzw$, where $w \in N_H(xyz, W_\textup{red}^*)$. 
We say that~$xyz$ generates the blue tight component~$B(xyz)$. It suffices to show that all $xyz \in \mathcal{T}$ generate the same blue tight component. 
First we show that triples that share two vertices generate the same blue tight component. 
Note that, for $xyz_1, xyz_2 \in \mathcal{T}$, we have $d_H(xyz_1, W_\textup{red}^*), d_H(xyz_2, W_\textup{red}^*) \geq \s{W_\textup{red}^*}-\eps n > \s{W_\textup{red}^*}/2$ and thus there exists $w \in N_H(xyz_1)\cap N_H(xyz_2) \cap W_\textup{red}^*$. 
Since the edges~$wxyz_1$ and~$wxyz_2$ are blue, it follows that $B(xyz_1) = B(xyz_2)$. 

Now let $x_1y_1z_1, x_2y_2z_2 \in \mathcal{T}$, where $x_1y_1, x_2y_2 \in G^{\red}$. Let $w_1 \in N_{\partial R}(x_1y_1) \cap N_{\partial R}(x_2y_2) \cap N_{G^{\red}}(x_1)\cap N_{G^{\red}}(x_2) \cap W_\textup{red}^*$ and $w_2 \in N_{\partial R}(x_1w_1) \cap N_{\partial R}(x_2w_1) \cap W_\textup{red}^*$.
It follows that $x_1y_1w_1, x_1w_1w_2, x_2w_1w_2, x_2y_2w_1 \in \mathcal{T}$. Hence $B(x_1y_1z_1) = B(x_1y_1w_1) = B(x_1w_1w_2) = B(x_2w_1w_2) = B(x_2y_2w_1) = B(x_2y_2z_2)$. Let~$B'$ be the unique blue tight component generated by all triples $xyz \in \mathcal{T}$.
\end{proofclaim}

The previous claim together with a greedy argument implies that there is a matching~$M_*^{B'}$ in~$B'[W^*]$ that covers all but~$\eta n$ of the vertices in~$W^*$. Thus we may assume that $\s{M_*^{\blue}}\geq \eta n/4$, otherwise $\s{V(M_*^{\red}\cup M_*^{B'})} \geq n - 3\eta n $, which is a contradiction to the initial assumption.
To complete the proof, we will show that in fact~$B' =B$, implying~$M_*^{\red}$ and $M_*^{\blue} \cup M_*^{B'}$ are tightly connected matchings, a contradiction to the initial assumption.

We now pick a special edge $e^* \in M_*^{\blue}$. 
Its special property that we desire is stated in~\cref{claim:edge restriction}.


\begin{claim}
\label{claim:e^*}
There exist an edge $e^* = v_1^*v_2^*v_3^*v_4^* \in M_*^{\blue}$ and distinct vertices $w_1, \dots, w_4,$ $w_1', \dots, w_4' \in W_\textup{red}^*$ such that, for each $j \in [4]$, 
\begin{enumerate}[label = \upshape (\alph*)]
\item all the red edges of~$G$ incident to~$v_j^*$ induce~$R$, or \label{itm:6.5-a}
\item $v_j^*w_j \in G^{\blue}$ and $v_j^*w_jw_j' \in \partial R \cap \partial B$. \label{itm:6.5-b}
\end{enumerate}
\end{claim}
\begin{proofclaim}
For each edge $e \in M_*^{\blue}$, let $v_1^e,v_2^e,v_3^e,v_4^e$ be an enumeration of its vertices.
It is easy to see that there exists $M_1^{\blue} \subseteq M_*^{\blue}$ with $\s{ M_1^{\blue}} = \s{ M_*^{\blue} }/16$ such that for each $j\in [4]$ we have that either
\begin{enumerate}[label = \upshape (\alph*$ '$)]
    \item for all $e \in M_1^{\blue}$, there is a red edge in~$G$ between~$v_j^e$ and~$W_\textup{red}^*$, or \label{1}
    \item for all $e \in M_1^{\blue}$, all edges in~$G$ between~$v_j^e$ and~$W_\textup{red}^*$ are blue.
\end{enumerate}
Let~$J_1$ be the set of $j \in [4]$ such that \ref{1} holds and $J_2 = [4]\setminus J_1$.
Since each vertex in~$W_\textup{red}^*$ is incident to a red edge of~$G$ that induces~$R$ and~$G$ is a blueprint for~$H$, we have that, for all $e \in M_1^{\blue}$ and all $j \in J_1$, all the red edges incident to~$v_j^e$ induce~$R$.
For every $j \in J_2$, we have that 
\begin{align*}
    \s{ G^{\blue} [\left\{v_j^e \colon e \in M_1^{\blue}\right\}, W_\textup{red}^*] } &\geq \s{ M_1^{\blue} } \br{ \s{ W_\textup{red}^*} - \alpha^{1/30}n } 
    \geq \br{1-\alpha^{1/31}} \s{ M_1^{\blue} } \s{ W_\textup{red}^*}.
\end{align*}
Thus there exists $w_j \in W_\textup{red}^*$ such that~$w_jv_j^e$ is blue for at least $\s{M_1^{\blue}}(1-\alpha^{1/32})$ of the vertices~$v_j^e$, with $e \in M_1^{\blue}$.
It is easy to see that we can choose the~$w_j$ to be distinct. Hence there exist distinct vertices $w_1,w_2,w_3,w_4 \in W_\textup{red}^*$ and $M_2^{\blue} \subseteq M_1^{\blue}$ with $\s{M_2^{\blue}} = \s{M_1^{\blue}}/2 \geq \eta n /128$ such that for all $j \in J_2$ and all $e \in M_2^{\blue}$ we have that $w_jv_j^e \in G^{\blue}$. 

For $j \in J_2$, let $V_j= \{v_j^e \colon e \in M_2^{\blue}\}$ and note that $d_{G}^{\blue}(w_j,V_j) = \s{M_2^{\blue}} \geq \eta n /128$ and $ d_{G}^{\red}(w_j, W_\textup{red}^*) \geq \eta n/2$. 
For each $j \in J_2$, we apply \cref{lem:shadow} with colours reversed and $w_j, V_j, \widetilde{W}_\textup{red}^*$ playing the roles of $T, S^{\blue}, S^{\red}$ where $\widetilde{W}_\textup{red}^*$ denotes~$W_\textup{red}^*$ with all previously chosen vertices removed. 
Thus, we find distinct $w_j' \in W_\textup{red}^*\setminus \{w_1,w_2,w_3,w_4\}$ and $M_3^{\blue} \subseteq M_2^{\blue}$ with $\s{M_3^{\blue}} = \s{M_2^{\blue}}/2$ such that, for all $j \in J_2$ and all $e \in M_3^{\blue}$, we have that $v_j^ew_j \in G^{\blue}$ and $v_j^ew_jw_j' \in \partial R \cap \partial B$.
We complete the proof by choosing $e^* = v_1^*v_2^*v_3^*v_4^* \in M_3^{\blue}$ and a distinct vertex $w_j' \in W_\textup{red}^*$ for each $j \in J_1$.
\end{proofclaim}

Let $W' = W_\textup{red}^* \setminus \{w_1, \dots, w_4, w_1', \dots, w_4'\}$.

\begin{claim}
\label{claim:edge restriction}
The $4$-graph $R[e^* \cup W']$ is empty and $B[e^* \cup W']$ does not contain two vertex-disjoint edges each containing an edge of~$G^{\blue}$.
In particular, there do not exist two vertex-disjoint edges~$f_1$ and~$f_2$ in $(R \cup B)[e^* \cup W']$ each containing an edge of~$G^{\blue}$. 
\end{claim}

\begin{proofclaim}
First suppose there exist two vertex-disjoint edges $f_1, f_2 \in B[e^* \cup W']$ each of which contains an edge of~$G^{\blue}$. 
By the maximality of~$\s{M_*}$, both~$f_1$ and~$f_2$ must intersect~$e^*$. 
For simplicity, we only consider the case that $e^* \setminus (f_1 \cup f_2) = \{v_1^*\}$ (the other cases can be proved similarly).
By \cref{claim:e^*}, we have that all red edges of~$G$ incident to~$v_1^*$ induce~$R$ or $v_1^*w_1 \in G^{\blue}$ and $v_1^*w_1w_1' \in \partial R \cap \partial B$.

First suppose that $v_1^*w_1 \in G^{\blue}$ and $v_1^*w_1w_1' \in \partial R \cap \partial B$. 
Let $w_1'' \in N_H(v_1^*w_1w_1',W^*\setminus (f_1\cup f_2))$ and $f_3 = v_1^*w_1w_1'w_1''$. 
Let $M' = (M_*\setminus \{e^*\})\cup \{f_1,f_2,f_3\}$. 
Note that $W(M') \subseteq W^*$. Since $ |W| \ge \eta n \ge 3 \gamma n $ and $|W^*_\textup{blue}| \le \gamma n $ by~\cref{claim:mathcalM}, we deduce that~$M'$ satisfies~\ref{itm:M4}.
Hence $M' \in \mathcal{M}'$ contradicting the maximality of~$\s{M_*}$.

Now assume that all the red edges of~$G$ incident to~$v_1^*$ induce~$R$. 
Let~$M$ be a matching in $R \cup B$ containing $(M_* \setminus \{e^*\}) \cup \{f_1, f_2\}$ satisfying~\ref{itm:M2} and~\ref{itm:M3}.
We now show that $M \in \mathcal{M}'$, which then contradicts the maximality of~$\s{M_*}$. 
Recall that $v_1^* \in e^* \in M_*^{\blue}$, so
\begin{align}
W \subseteq (W^*\setminus (f_1 \cup f_2)) \cup \{v_1^*\} \text{ and }
V(M^{\blue}) \cup W \subseteq V(M_*^{\blue}) \cup W^*.
\label{eqn:WW_*}
\end{align}
Together with our assumption on~$v_1^*$,~$M$ satisfies~\ref{itm:M1}.
Hence $M \in \mathcal{M}$. 
For all $w \in W \cap W_\textup{red}^*$, 
\begin{align*}
	d_{G[W]}^{\blue}(w) 
	\overset{\eqref{eqn:WW_*}}{\leq} d_{G[W^*]}^{\blue}(w) +|v^*_1| 
	\overset{\text{\cref{claim:mathcalM}\ref{itm:Ma}}}{\le} 7 \sqrt{\eps}n + 1 
	\leq 8\sqrt{\eps}n.
\end{align*}
and a similar inequality holds for all $w \in W \cap W_\textup{blue}^*$.
This implies that $W_\textup{blue} \subseteq W^*_\textup{blue} \cup \{v^*\}$.
Since $ |W| \ge \eta n \ge 3 \gamma n $ and $|W^*_\textup{blue}| \le \gamma n $ by~\cref{claim:mathcalM}, we deduce that~$M$ satisfies~\ref{itm:M4}.
Hence, $M \in \mathcal{M'}$ as required, a contradiction. 

Therefore, $B[e^* \cup W']$ does not contain two vertex-disjoint edges each of which contains an edge of~$G^{\blue}$.
If $R[e^* \cup W']$ contains an edge~$f$, then a similar argument holds with~$f$ replacing $\{f_1, f_2\}$. 
Note that if $\s{M} = \s{M_*}$, then we obtain a contradiction by showing that $\s{M_*^{\red}} < \s{M^{\red}}$. 
\end{proofclaim}

Since $e^* \in M_*^{\blue} \subseteq (G^{\blue})^+$, we may assume without loss of generality that $v_1^*v_2^* \in G^{\blue}$. The following claim shows that one of the vertices~$v_1^*$ and~$v_2^*$ has small blue degree in~$G$ to~$W'$ (and thus it has large red degree to~$W'$). 
\begin{claim}
\label{claim:at most one}
We have $d_{G}^{\blue}(v_1^*, W') \leq 3\gamma n$ or $d_{G}^{\blue}(v_2^*,W') \leq 3\gamma n$.
\end{claim}
\begin{proofclaim}
Suppose to the contrary that we have $d_{G}^{\blue}(v_1^*, W'), d_{G}^{\blue}(v_2^*,W') > 3\gamma n$.
By \cref{claim:edge restriction}, it suffices to show that we can find two vertex-disjoint edges~$f_1$ and~$f_2$ in $(R \cup B)[e^* \cup W']$ each containing an edge of~$G^{\blue}$. 
It is easy to see that we can greedily choose vertices $x \in N_{G}^{\blue}(v_1^*, W')$, $x' \in N_{G}^{\red}(x,W') \cap N_{\partial B}(v_1^*x, W')$ and $x'' \in N_{\partial R}(xx', W') \cap N_H(v_1^*xx',W')$. Set $f_1 = v_1^*xx'x''$. By our construction, $v_1^*xx' \in \partial B$ and $xx'x'' \in \partial R$ implying $f_1 \in (R\cup B)[e^* \cup W']$. Similarly there exists an edge $f_2 = v_2^*yy'y'' \in (R \cup B)[e^* \cup W']$ disjoint from~$f_1$ with $y,y',y'' \in W'$.
\end{proofclaim}

Without loss of generality assume $d_{G}^{\blue}(v_1^*,W') \leq 3 \gamma n$ and so $d_{G}^{\red}(v_1^*, W') \geq \s{W'} - \alpha^{1/31}n$.
Let $w \in N_{\partial B}(v_1^*v_2^*) \cap N_{G}^{\red}(v_1^*) \cap W',$ 
$w' \in N_{G}^{\red}(w) \cap N_{\partial R}(v_1^*w) \cap N_H(v_1^*v_2^*w) \cap W'$ and $w'' \in N_H(v_1^*ww',W')$.
(We can find these vertices greedily one by one.)
By \cref{claim:B'}, we may further assume that $ww'w'' \in \partial B'$. 
By construction, we have that $v_1^*ww' \in \partial R$ and thus \cref{claim:edge restriction} implies that both~$v_1^*v_2^*ww'$ and~$v_1^*ww'w''$ are blue. Since $v_1^*v_2^*w \in \partial B$, we deduce that $v_1^*v_2^*ww', v_1^*ww'w'' \in B$ and so $ww'w'' \in \partial B$ implying that $\partial B \cap \partial B' \neq \varnothing$. Therefore~$B = B'$ as required.
\end{proof}

\section{Monochromatic connected matchings in \texorpdfstring{$K_n^{(5)}$}{Kn(5)}}
\label{section:matchings Kn5}

The aim of this section is to prove the following lemma which states that every $2$-edge-coloured dense $5$-graph can be almost partitioned into four monochromatic tightly connected matchings.

\begin{lem}
\label{lem:5-matchings}
Let $1/n \ll \eps \ll \alpha \ll \eta < 1$.
Let~$H$ be a $2$-edge-coloured $(1-\eps,\alpha)$-dense $5$-graph on~$n$ vertices. Then~$H$ contains four vertex-disjoint monochromatic tightly connected matchings such that their union covers all but at most~$3\eta n$ of the vertices of~$H$.
\end{lem}

Note that this implies $\mu_5^4(1, \eps, n) \geq (1- 3\eta)n/5$ for $1/n \ll \eps \ll \eta < 1$. Hence $\mu_5^4(1) \geq 1/5$.
Together with \cref{cor:matchings_to_cycles}, \cref{lem:5-matchings} implies \cref{thm:2}.

We use the following notation throughout this section.
Let~$H$ be a $2$-edge-coloured $5$-graph and let~$G$ be a blueprint for~$H$. Given a red tight component $R \subseteq H$, we write~$R^3$ for the edges of~$G$ that induce~$R$. We use analogous notation for blue tight components.

Let~$H$ be a $2$-edge-coloured dense $5$-graph.
We first apply~\cref{lem:generalblueprint} to~$H$ to get a blueprint~$G$ for~$H$. 
Since~$G$ is $2$-edge-coloured dense $3$-graph, we can apply \cref{lem:generalblueprint} again to~$G$ to obtain a blueprint for~$G$, which is a $2$-coloured $1$-graph. 
The following lemma summarises the structural information about~$H$ that we obtain in this way.

\begin{lem}
\label{lem:5-blueprint}
Let $1/n \ll \eps \ll \alpha \ll 1$. 
Let~$H$ be a $2$-edge-coloured $(1-\eps, \alpha)$-dense $5$-graph on~$n$ vertices.
Then there exists a $3$-graph~$G$ with $V(G) = V(H)$, two disjoint subsets~$V^{\red}$ and~$V^{\blue}$ of~$V(H)$, a red tight component $R \subseteq H$ and a blue tight component $B\subseteq H$ such that the following properties hold.
\begin{enumerate}[label = \upshape (\roman*)]
    \item$G$ is a $(1-\alpha^{1/37}, \alpha^{1/37})$-dense $3\sqrt{\eps}$-blueprint for~$H$.
    \item $\s{V(H) \setminus (V^{\red} \cup V^{\blue})} \leq \alpha^{1/75}n$.
    \item $d_{\partial R^3}(v) \geq (1-\alpha^{1/75})n$ for all $v \in V^{\red}$.
    \item $d_{\partial B^3}(v) \geq (1-\alpha^{1/75})n$ for all $v \in V^{\blue}$.
\end{enumerate}
\end{lem}
\begin{proof}
By \cref{lem:generalblueprint}, there exists a $(1-\alpha^{1/37}, \alpha^{1/37})$-dense $3\sqrt{\eps}$-blueprint~$G$ for~$H$ with $V(G) = V(H)$. 
We apply \cref{lem:generalblueprint} to~$G$ and obtain a $\alpha^{1/75}$-blueprint~$J$ for~$G$ with $\s{J} \geq (1-\alpha^{1/75})n$. 
Note that, as a blueprint for a $3$-graph,~$J$ is a $1$-graph. Hence each edge of~$J$ contains precisely one vertex. By the definition of a blueprint all the red edges of~$J$ induce the same red tight component~$R_G$ of~$G$. Let $V^{\red} = \bigcup J^{\red}$. Since~$R_G$ is a red tight component of~$G$ all its edges induce the same red tight component~$R$ of~$H$. Define~$V^{\blue}$ and~$B$ analogously.
\end{proof}

Two edges~$f$ and~$f'$ in~$H$ are \emph{loosely connected} if there exists a sequence of edges $e_1,\dots,e_t$ such that~$e_1 = f$,~$e_t=f'$ and $\s{e_i \cap e_{i+1}} \geq 1$ for all $i \in [t-1]$.
A subgraph~$H'$ of~$H$ is \emph{loosely connected} if every pair of edges in~$H'$ is loosely connected. A maximal loosely connected subgraph of~$H$ is called a \emph{loose component} of~$H$.

We now prove \cref{lem:5-matchings}.
The proof works by first finding a maximal matching in $R \cup B$, where~$R$ and~$B$ are the components given by \cref{lem:5-blueprint}, and then finding maximal connected matchings in the remaining vertices. 

\begin{proof}[Proof of \cref{lem:5-matchings}]
Assume, for a contradiction, that such matchings do not exist. We call this the initial assumption. Apply \cref{lem:5-blueprint} and obtain $V^{\red}, V^{\blue}, G, R^3, R, B^3, B$ and let $V^* = V^{\red} \cup V^{\blue}$. Since there are only few vertices in $V(H) \setminus V^*$ we ignore these vertices from the start and construct our matchings in~$H[V^*]$.

We begin by choosing a matching $M \subseteq (R \cup B)[V^*]$ of maximum size. Let $U = V^*\setminus V(M)$. Note that we have $R[U] = B[U] =\varnothing$ and $\s{U} \geq \eta n$ by the initial assumption. Let $U^{\red} = U \cap V^{\red}$ and $U^{\blue} = U \cap V^{\blue}$.
The following claim shows that if~$U^{\red}$ and~$U^{\blue}$ are both large, then~$G[U]$ must contain many edges in~$R^3$ or many edges in~$B^3$.

\begin{claim}
\label{claim:structure}
If $\s{U^{\red}}, \s{U^{\blue}} \geq \alpha^{1/309}n$, then $\max\{\s{R^3[U]},\s{B^3[U]}\} \geq \frac{1}{2}\s{U^{\red}}\s{U^{\blue}}\s{U} - 3\alpha^{1/155}n^3$.
\end{claim}
\begin{proofclaim}
Define a bipartite graph~$K_0$ with vertex classes~$U^{\red}$ and~$U^{\blue}$ such that $x \in U^{\red}$ and $y \in U^{\blue}$ are joined by an edge if and only if $xy \in \partial R^3 \cap \partial B^3$. Recall that $d_{\partial R^3}(x) \geq (1-\alpha^{1/75})n$ and $d_{\partial B^3}(y) \geq (1-\alpha^{1/75})n$ for all $x \in U^{\red}$ and $y \in U^{\blue}$.
Hence
\[
\s{K_0} \geq \s{U^{\blue}}\s{U^{\red}} - \alpha^{1/75}n^2.
\]
Since~$G$ is $(1-\alpha^{1/37}, \alpha^{1/37})$-dense, we have $d_G(xy, U) \geq \s{U} - \alpha^{1/37}n$ for $xy \in K_0$.
We now colour the edges of~$K_0$ such that $xy \in K_0$ is red if $d_{R^3}(xy, U) \geq \s{U} - 2\alpha^{1/76}n$ and blue if $d_{B^3}(xy, U) \geq \s{U} - 2\alpha^{1/76}n$. 
Since $K_0 \subseteq \partial R^3 \cap \partial B^3$, if $xyz \in G$ with $xy \in K_0$, then $xyz \in R^3 \cup B^3$.
Hence it suffices to show that almost all edges of~$K_0$ are of the same colour. Indeed, if we have that at least $\s{U^{\red}}\s{U^{\blue}} - 3 \alpha^{1/154}n^2$ edges of~$K_0$ are red, then we have 
\[
    \s{R^3[U]} \geq \frac{1}{2}(\s{U^{\red}}\s{U^{\blue}} - 3\alpha^{1/154}n^2) (\s{U} - 2 \alpha^{1/76}n) \geq \frac{1}{2} \s{U^{\red}}\s{U^{\blue}}\s{U} - 3\alpha^{1/155}n^3.
\]

We show that each edge $xy \in K_0$ is coloured either red or blue. It suffices to show that either $d_{R^3}(xy, U) < \alpha^{1/76}n$ or $d_{B^3}(xy, U) < \alpha^{1/76}n$. Indeed if $d_{R^3}(xy, U), d_{B^3}(xy, U) \geq \alpha^{1/76}n$, then by \cref{lem:shadow}, there exists $u, u' \in U$ such that $xyu \in R^3$, $xyu' \in B^3$ and $xyuu' \in \partial R \cap \partial B$. For any $u'' \in N_H(xyuu', U)$, we would have $xyuu'u'' \in R[U] \cup B[U]$, a contradiction to the maximality of~$M$.
Moreover, by \cref{lem:vertexdeg}, we have that $\min \{ d_{K_0}^{\red}(u), d_{K_0}^{\blue}(u)\} \leq \alpha^{1/76}n$ for all $u \in U$.

Let~$K_1$ be the graph obtained from~$K_0$ by, for each~$u \in U$, deleting all red edges incident to~$u$ if $d_K^{\red}(u) \leq \alpha^{1/76}n$ and all blue edges incident to~$u$ if $d_K^{\blue}(u) \leq \alpha^{1/76}n$. Note that $\s{K_1} \geq \s{U^{\red}}\s{U^{\blue}} - \alpha^{1/77}n^2$ and that, in~$K_1$, each vertex is incident to only edges of one colour. It is not too hard to see that by deleting at most $2\alpha^{1/154}n^2$ additional edges, we can obtain a subgraph~$K_2$ of~$K_1$ for which each vertex has degree~$0$ or large degree. More precisely, for all $u \in U^{\red}$,
\[ 
d_{K_2}(u) \geq \s{U^{\blue}} - 3\alpha^{1/308}n \text{ or } d_{K_2}(u) = 0
\] 
and, for all $u \in U^{\blue}$,
\[ 
d_{K_2}(u) \geq \s{U^{\red}} - 3\alpha^{1/308}n \text{ or } d_{K_2}(u) = 0.
\] 
Since each vertex is incident to only edges of one colour and any two vertices in~$U^{\red}$ that have non-zero degree have a common neighbour this implies that all edges in~$K_2$ are of the same colour. Since $\s{K_2} \geq \s{U^{\red}}\s{U^{\blue}} - 3\alpha^{1/154}n^2$, this concludes the proof.
\end{proofclaim}

The following claim shows that there is a red tight component~$R_*$ and a blue tight component~$B_*$ of~$H$ such that almost all the edges in~$G[U]$ induce one of these components.

\begin{claim}
\label{claim:two_components}
Let $\gamma = \alpha^{1/1110}$.
There exists a red tight component~$R_*$ and a blue tight component~$B_*$ of~$H$ such that 
\begin{enumerate}[label = \upshape (\roman*)]
    \item \label{two_components_i} $\s{R_*^3[U]} \geq \s{G^{\red}[U]} - 8 \gamma^{1/5}n^3$ and
    $\s{B_*^3[U]} \geq \s{G^{\blue}[U]} - 8 \gamma^{1/5}n^3$,
    \item \label{two_components_ii} $\s{(R_*^3 \cup B_*^3)[U]} \geq (1-\gamma^{1/6})\binom{\s{U}}{3}$ and
    \item \label{two_components_iii} $R_* = R$ or $B_* = B$.
\end{enumerate}
\end{claim}
\begin{proofclaim}
First we show that, for each $u \in U$, there exists $J_u \subseteq G_u[U]$, where~$G_u$ is the link graph of~$G$ at~$u$, such that $\s{J_u} \geq \s{G_u[U]} - \alpha^{1/14}n^2$ and $R(e \cup u) = R(e' \cup u)$ for $e, e' \in J_u^{\red}$ and $B(e \cup u) = B(e' \cup u)$ for $e, e' \in J_u^{\blue}$. 

To show this fix $u \in U$. Without loss of generality assume that $u \in U^{\red}$. By \cref{lem:5-blueprint}, $d_{\partial R^3}(u,U) \geq \s{U} - \alpha^{1/75}n$. 
Let $U_* = N_{\partial R^3}(u,U)$. Clearly, $\s{U_*} \geq \eta n/2$ and $G_u^{\red}[U_*] \subseteq R^3_u$.
Moreover, for all $x \in U_*$, we have $d_G(ux) > 0 $ and thus, since~$G$ is $(1-\alpha^{1/37}, \alpha^{1/37})$-dense, $d_G(ux) \geq (1-\alpha^{1/37})n$. It follows that $\delta(G_u[U_*]) \geq \s{U_*} - \alpha^{1/37}n$.
Thus by applying \cref{lem:reducing_components} with $R^3, u, U_*, \alpha^{1/37}$ playing the roles of $R_G, S, U, \delta$, there exists $J_u \subseteq G_u[U_*] \subseteq G_u[U]$ such that 
\[
\s{J_u} \geq \s{G_u[U_*]} - 7\alpha^{1/148}n^2 \geq \s{G_u[U]} - \alpha^{1/75}n^2 - 7\alpha^{1/148}n^2 \geq \s{G_u[U]} - \alpha^{1/149}n^2
\]
and 
$H(u \cup e) = H(u \cup e')$ for $e, e' \in J_u$ of the same colour.

Now consider the auxiliary multi-$3$-graph $D = \bigcup_{u \in U} \{ e \cup u \colon e \in J_u\}$. 
Note that
\[
\s{D} = \sum_{u \in U} \s{J_u} \geq \sum_{u \in U}\left(\s{G_u[U]} - \alpha^{1/149}n^2\right) \geq 3\s{G[U]} - \alpha^{1/149}n^3.
\]
Let~$F$ be the subgraph of~$G[U]$ for which $e \in F$ if and only if~$e$ is an edge of multiplicity~$3$ in~$D$.
Since~$G$ is $(1-\alpha^{1/37}, \alpha^{1/37})$-dense, \cref{prop:edges_in_dense} implies that $\s{G} \geq (1 - 2 \alpha^{1/37}) \binom{n}{3}$. Hence
\begin{align*}
    \s{G[U]} &\geq \binom{\s{U}}{3} - 2\alpha^{1/37}\binom{n}{3} 
    \geq \binom{\s{U}}{3} - 2\alpha^{1/37} \binom{\s{U}/\eta}{3} \\
    &\geq \binom{\s{U}}{3} - \frac{4\alpha^{1/37}}{\eta^3}\binom{\s{U}}{3}
    \geq (1-\alpha^{1/38})\binom{\s{U}}{3}.
\end{align*}
Therefore $\s{F} \geq \s{G[U]} - \alpha^{1/149}n^3 \geq (1-\alpha^{1/150})\binom{\s{U}}{3}$. Recall that $\gamma = \alpha^{1/1110}$.
By \cref{prop:dense,prop:edges_in_dense}, there exists a $(1-\gamma^{1/5}, \gamma^{1/5})$-dense subgraph $\tF \subseteq F$ with $V(\tF) = V(F) =U$ and, by \cref{prop:edges_in_dense}, $\s{\tF} \geq (1-2\gamma^{1/5}) \binom{\s{U}}{3}$. Hence $\s{\tF^{\red}} \geq \s{G^{\red}[U]} - 2\gamma^{1/5}n^3$.
Let $S^{\red} = \{x \in U \colon d_{\tF^{\red}}(x) \geq 6 \gamma^{1/5} n^2\}$.
Let~$F_0^{\red}$ be the subgraph of~$\tF^{\red}$ consisting of all edges that contain a vertex in~$S^{\red}$. Note that $\s{F_0^{\red}} \geq \s{\tF^{\red}} - 6 \gamma^{1/5}n^3 \geq \s{G^{\red}[U]} - 8\gamma^{1/5}n^3$.

We claim that all the edges in~$F_0^{\red}$ induce the same red tight component~$R_*$ in~$H$.
Let $e, e' \in \tF^{\red}$ with $u \in e \cap e'$. Note that $e \setminus u, e' \setminus u \in J_u^{\red}$ and so $R(e) = R(e')$. Hence edges in the same loose component of~$\tF^{\red}$ induce the same red tight component in~$H$. In particular, since $F_0^{\red} \subseteq \tF^{\red}$, for $u \in S^{\red}$, all the edges in $N_{F_0^{\red}}(u)$ induce the same red tight component~$R(u)$ of~$H$.

Let $u, v \in S^{\red}$. We want to show that $R(u) = R(v)$. We may assume that~$u$ and~$v$ are in distinct loose components~$L$ and~$L'$ of~$\tF^{\red}$, respectively. In particular, any edge of~$\tF$ that intersects both~$V(L)$ and~$V(L')$ is in~$\tF^{\blue}$.
If $u, v \in V^{\red}$, then $d_{\partial R^3}(u), d_{\partial R^3}(v) \geq (1-\alpha^{1/75})n$ implying $R(u) = R =R(v)$. 
Thus we may assume that one of~$u$ and~$v$ is in~$V^{\blue}$, say $v \in V^{\blue}$.
Let $\Gamma_L(u) = \{ u' \in V(L) \colon d_L(uu') \geq \gamma^{1/5}n\}$ and $\Gamma_{L'}(v) = \{ v' \in V(L') \colon d_{L'}(vv') \geq \gamma^{1/5}n\}$. 
It is easy to see that $\s{\Gamma_L(u)}, \s{\Gamma_{L'}(v)} \geq 5 \gamma^{1/5}n$. Let~$D'$ be the bipartite directed graph with parts $\Gamma_L(u)$ and $\Gamma_{L'}(v)$ such that, for $u' \in \Gamma_L(u)$,
\begin{align*}
N_{D'}^+(u') = \{v' \in \Gamma_{L'}(v) \colon &uu'v' \in \tF^{\blue} \text{ and } uu'u''v' \in \partial R(uu'u'') \cap \partial B(uu'v') \\ &\text{and } uu'u''vv' \in H \text{ for some } u'' \in N_L(uu')\}, 
\end{align*}
and, for $v' \in \Gamma_{L'}(v)$,
\begin{align*}
N_{D'}^+(v') = \{u' \in \Gamma_{L}(u) \colon &vv'u' \in \tF^{\blue} \text{ and } u'v \in \partial B^3 \text{ and } vv'v''u' \in \partial B \cap \partial R(vv'v'') \\ &\text{and } vv'v''uu' \in H \text{ for some } v'' \in N_{L'}(vv')\}.
\end{align*}
By \cref{lem:shadow}, the fact that $\widetilde{F}$ is $(1-\gamma^{1/5}, \gamma^{1/5})$-dense and the fact that~$H$ is $(1-\eps, \alpha)$-dense, we have, for $u' \in \Gamma_L(u)$,
\[
d_{D'}^+(u') \geq \s{\Gamma_{L'}(v)} - \gamma^{1/5}n - \eps^{1/4}n - \eps n > \s{\Gamma_{L'}(v)}/2.
\]
Similarly, also using the fact that $d_{\partial B^3}(v) \geq (1-\alpha^{1/75})n$, we have, for $v' \in \Gamma_{L'}(v)$,
\[
d_{D'}^+(v') \geq \s{\Gamma_{L}(u)} - \gamma^{1/5}n -\alpha^{1/75}n - \eps^{1/4}n - \eps n > \s{\Gamma_{L}(u)}/2.
\]
It follows that~$D'$ contains a double edge~$u'v'$, where $u' \in \Gamma_{L}(u)$ and $v' \in \Gamma_{L'}(v)$. Let $u'' \in N_{L}(uu')$ and $v'' \in N_{L'}(vv')$ be the vertices that are guaranteed to exist by the definition of~$D'$. Since $u'v \in \partial B^3$, we have that $vv'u \in B^3$ and thus also $uu'v' \in B^3$. As $B[U] = \noth$, we have $vv'v''uu', uu'u''vv' \in H^{\red}$ and thus $R(uu'u'') = R(vv'v'')$.
Hence $R(u) = R(v)$.
We define~$F_0^{\blue}$ and~$B_*$ in an analogous way. This proves \ref{two_components_i}.

Note that \ref{two_components_ii} follows from \ref{two_components_i} using the facts $\s{U} \geq \eta n$ and $\s{G[U]} \geq (1- \alpha^{1/38})\binom{\s{U}}{3}$, which were noted earlier in this proof. 

We will now prove \ref{two_components_iii}. We distinguish between two cases.
\begin{enumerate}[label=\textbf{Case \arabic*:},wide, labelwidth=0pt,labelindent=0pt,parsep=0pt]
\item \boldmath $\s{U^{\red}}, \s{U^{\blue}} \geq \gamma^{1/13}n.$ \unboldmath \\
By \cref{claim:structure}, we have $\max\{\s{R^3[U]}, \s{B^3[U]}\} \geq \frac{1}{2} \s{U^{\red}}\s{U^{\blue}}\s{U} - 3\alpha^{1/155}n^3$. Since $\frac{1}{2} \s{U^{\red}}\s{U^{\blue}}\s{U} - 3\alpha^{1/155}n^3 \geq \frac{1}{2}\gamma^{2/13}\eta n^3 - 3 \alpha^{1/155}n^3 \geq 2 \gamma^{1/6}n^3$, we have $R_*^3 \cap R^3 \neq \varnothing$ or $B_*^3 \cap B^3 \neq \varnothing$ and thus $R_* = R$ or $B_* = B$.

\item \boldmath $\s{U^{\blue}} \leq \gamma^{1/13}n$ \textbf{or} $\s{U^{\red}} \leq \gamma^{1/13}n.$ \unboldmath \\
Say $\s{U^{\blue}} \leq \gamma^{1/13}n$. Then $\s{U^{\red}} = \s{U} - \s{U^{\blue}} \geq \s{U} - \gamma^{1/13}n$. Let $Q^3 = \{ T \in \binom{U}{3} \colon \binom{T}{2} \cap \partial R^3 \neq \varnothing\}$. Since $d_{\partial R}(u, U) \geq \s{U} - \alpha^{1/75}n$ for $u \in U^{\red}$, there can be at most $\s{U^{\red}}\alpha^{2/75}n^2$ triples that intersect~$U^{\red}$ and are not in~$Q^3$. Hence
\begin{align*}
    \s{Q^3} &\geq \binom{\s{U}}{3} - \s{U^{\blue}}^3 - \s{U^{\red}}\alpha^{2/75}n^2 \\
    &\geq \binom{\s{U}}{3} - \gamma^{3/13}n^3 - \alpha^{2/75}n^3 
    \geq \binom{\s{U}}{3} - 2\gamma^{1/5}n^3.
\end{align*}
Note that $\s{R^3[U]} \geq \s{Q \cap G^{\red}[U]} \geq \s{G^{\red}[U]} - 2\gamma^{1/5}n^3$.
Therefore, we have $R_* = R$. \qedhere
\end{enumerate}
\end{proofclaim}

We define $R_\diamond = R \cup R_*$ and $B_\diamond = B \cup B_*$.
Note that, by \cref{claim:two_components}\ref{two_components_iii}, $R_\diamond \cup B_\diamond$ is the union of at most three monochromatic tight components.
Let~$M_\diamond$ be a maximal matching in $(R_\diamond \cup B_\diamond)[V^*]$ containing~$M$. Let $W = V^*\setminus V(M_\diamond)$. Since $M \subseteq M_\diamond$, we have $W \subseteq U$. By the initial assumption, we have $\s{W} \geq \eta n$. Note that $(R_* \cup B_*)[W] = \varnothing$ and, since $W \subseteq U$, $(R_*\cup B_*)[W] \geq \binom{\s{W}}{3}- \gamma^{1/6}n^3$.
The following claim shows that almost all the edges in~$G[W]$ are of the same colour.

\begin{claim}
\label{claim:one_component}
We have $\s{R_*^3[W]} \geq \binom{\s{W}}{3} - \gamma^{1/9}n^3$ or $\s{B_*^3[W]} \geq \binom{\s{W}}{3} - \gamma^{1/9}n^3$.
\end{claim}
\begin{proofclaim}
Let $G_* = R_*^3 \cup B_*^3$. 
We define
\begin{align*}
    W_\textup{red} = \{u \in W &\colon d_{G_*}(u,W) \geq 2\alpha n^2\text{ and } d_{B_*^3}(u,W) < \alpha n^2\}, \\
    W_\textup{blue} = \{u \in W &\colon d_{G_*}(u,W) \geq 2\alpha n^2\text{ and } d_{R_*^3}(u,W) < \alpha n^2\}, \\
    W_0 = \{u \in W &\colon d_{G_*}(u, W) < 2 \alpha n^2\}.
\end{align*}
Since $(R_* \cup B_*)[W] = \varnothing$, by \cref{lem:vertexdeg}, $W_\textup{red}, W_\textup{blue}$ and~$W_0$ partition~$W$.
Let~$J$ be the subgraph of~$G_*[W]$ obtained by deleting all red edges containing a vertex in $W_\textup{blue} \cup W_0$ and all blue edges containing a vertex in $W_\textup{red} \cup W_0$.
Note that $\s{J} \geq \s{G_*[W]} - 2\alpha n^3 \geq (1- \gamma^{1/7})\binom{\s{W}}{3}$ and $J \subseteq \binom{W_\textup{red}}{3} \dot\cup \binom{W_\textup{blue}}{3}$.
Hence 
\begin{align}
\label{ub}
     (1-\gamma^{1/7})\binom{\s{W}}{3} \leq \binom{\s{W_\textup{red}}}{3} + \binom{\s{W_\textup{blue}}}{3}.
\end{align}
Suppose that $\s{W_\textup{red}}, \s{W_\textup{blue}} \leq (1-\alpha^{1/8}) \s{W}$.
By (\ref{ub}), we may assume without loss of generality assume that $\s{W_\textup{red}} \geq \s{W}/2$. Noting that $x \mapsto x^3 + (\s{W} -x)^3 $ is an increasing function for $x \geq \s{W}/2$ we have 
\begin{align*}
    \binom{\s{W_\textup{red}}}{3} + \binom{\s{W_\textup{blue}}}{3} &\leq \frac{1}{6}\br{\s{W_\textup{red}}^3 + \s{W_\textup{blue}}^3} 
    \leq \frac{1}{6}\br{\s{W_\textup{red}}^3 + \br{\s{W} - \s{W_\textup{red}}}^3} \\
    &\leq ((1-\gamma^{1/8})^3 + \gamma^{3/8})\frac{\s{W}^3}{6} 
    < (1-\gamma^{1/7})\binom{\s{W}}{3},
\end{align*}
a contradiction to (\ref{ub}).

Hence at least one of $W_\textup{red}$ and $W_\textup{blue}$ has size at least $(1-\gamma^{1/8})\s{W}$. Without loss of generality assume $\s{W_\textup{red}} \geq (1 - \gamma^{1/8})\s{W}$. Note that any edge of~$J$ contained in $W_\textup{red}$ is in~$R_*^3$, hence 
\begin{align*}
    \s{R^3_*[W]} \geq \s{J} - \s{W \setminus W^{\red}}n^2 \geq \binom{\s{W}}{3} - \gamma^{1/9}n^3. 
\end{align*}
This proves the claim.
\end{proofclaim}

 Now assume without loss of generality that $\s{R_*^3[W]} \geq \binom{\s{W}}{3} - \gamma^{1/9}n^3$. 
Note that almost all edges in~$H[W]$ are blue (otherwise there would have to be an edge in~$R_*[W]$, which would contradict the maximality of~$M$). More precisely, we have
\[
\s{H^{\blue}[W]} \geq \frac{3!}{5!} \s{R_*[W]}(\s{W} - 3\sqrt{\eps}n)(\s{W} - \eps n) \geq (1-\gamma^{1/10})\binom{\s{W}}{5}.
\]
By \cref{prop:edges_in_dense,prop:dense}, there exists a $(1-\gamma^{1/1010}, \gamma^{1/1010})$-dense tightly connected subgraph $\widetilde{H}^{\blue}$ of $H^{\blue}[W]$ with $V(\widetilde{H}^{\blue}) = W$ and $\s{\widetilde{H}^{\blue}} \geq (1 - 2\gamma^{1/1010})\binom{\s{W}}{5}$. By an easy greedy argument, there exists a matching~$M'$ in $\widetilde{H}^{\blue}$ that covers all but at most $\eta n$ of the vertices in~$W$. The matching $M' \cup M_\diamond$ covers all but at most $3\eta n$ of the vertices of~$H$. This contradicts the initial assumption. \qedhere
\end{proof}

\section{Concluding Remarks} \label{sec:concluding}

For $k \ge 3$, let~$f(k)$ be the minimum integer~$m$ such that, for all large $2$-edge-coloured complete $k$-graphs, there exists~$m$ vertex-disjoint monochromatic tight cycles covering almost all vertices. 
Note that~$f(k)$ is well defined by~\cite{Bustamante2020} but the bound is very large. 
It is easy to see that $f(k) \geq 2$ for all $k \geq 3$. Indeed, consider the $k$-graph~$H = K^{(k)}(A,B)$ given in \cref{ex:blueprint} with $\s{A} = \frac{3k-1}{3k}n$. Note that~$H[A]$ is a red tight component. Moreover, note that any tight cycle contained in a monochromatic tight component other than~$H[A]$ covers at most about a third of the vertices of~$H$ and any tight cycle in~$H[A]$ leaves all $\frac{n}{3k}$ vertices in~$B$ uncovered. Hence no monochromatic tight cycle covers almost all vertices in~$H$.
We have $f(3) =2$ by~\cite{Bustamante2017}.
\cref{thm:1,thm:2} imply $f(4) = 2$ and $f(5) \leq 4$, respectively.
In general, we believe that $f(k) = 2$ for all~$k$. 
However, we believe that new ideas may be needed as indicated by again considering the $k$-graph~$H = K^{(k)}(A,B)$ with $\s{A} = \frac{3k-1}{3k}n$ (as above). 
If~$H$ contains two vertex-disjoint monochromatic tight cycles of distinct colour covering almost all vertices, then one of the two cycles must lie entirely in the red tight component~$H[A]$. 
However, this tight component is not induced by any edge in the blueprint of $H$ (which is $K^{(k-2)}(A,B)$ with colours swapped).
Thus we ask the weaker question of whether one can bound~$f(k)$ by some suitable function of~$k$.

\section*{Acknowledgements}
We thank Richard Lang, Nicol\'as Sanhueza-Matamala, and two anonymous referees for their helpful comments.


\begin{thebibliography}{10}

\bibitem{Allen2008}
P.~Allen.
\newblock Covering two-edge-coloured complete graphs with two disjoint
  monochromatic cycles.
\newblock {\em Combin. Probab. Comput.}, 17(4):471–486, 2008.

\bibitem{Allen2017}
P.~Allen, J.~B\"{o}ttcher, O.~Cooley, and R.~Mycroft.
\newblock Tight cycles and regular slices in dense hypergraphs.
\newblock {\em J. Combin. Theory Ser. A}, 149:30--100, 2017.

\bibitem{Bessy2010}
S.~Bessy and S.~Thomassé.
\newblock Partitioning a graph into a cycle and an anticycle, a proof of
  lehel's conjecture.
\newblock {\em J. Combin. Theory Ser. B}, 100(2):176--180, 2010.

\bibitem{Bustamante2020}
S.~Bustamante, J.~Corsten, N.~Frankl, A.~Pokrovskiy, and J.~Skokan.
\newblock Partitioning edge-colored hypergraphs into few monochromatic tight
  cycles.
\newblock {\em SIAM J. Discrete Math.}, 34(2):1460--1471, 2020.

\bibitem{Bustamante2017}
S.~Bustamante, H.~H\`an, and M.~Stein.
\newblock Almost partitioning 2-colored complete 3-uniform hypergraphs into two
  monochromatic tight or loose cycles.
\newblock {\em J. Graph Theory}, 91(1):5--15, 2019.

\bibitem{Erdos1991}
P.~Erdős, A.~Gyárfás, and L.~Pyber.
\newblock Vertex coverings by monochromatic cycles and trees.
\newblock {\em J. Combin. Theory Ser. B}, 51(1):90--95, 1991.

\bibitem{Garbe2019}
F.~Garbe, R.~Mycroft, R.~Lang, A.~Lo, and N.~Sanhueza-Matamala.
\newblock Partitioning \(2\)-coloured complete \(3\)-uniform hypergraphs into
  two monochromatic tight cycles.
\newblock {\em In preparation}.

\bibitem{Gyarfas2012}
A.~Gy\'{a}rf\'{a}s and G.~N. S\'{a}rk\"{o}zy.
\newblock Star versus two stripes {R}amsey numbers and a conjecture of
  {S}chelp.
\newblock {\em Combin. Probab. Comput.}, 21(1-2):179--186, 2012.

\bibitem{gyarfas2013}
A.~Gy{\'a}rf{\'a}s and G.~N. S{\'a}rk{\"o}zy.
\newblock Monochromatic path and cycle partitions in hypergraphs.
\newblock {\em Electron. J. Combin.}, 20(1):18, 2013.

\bibitem{Gyarfas2016}
A.~Gyárfás.
\newblock Vertex covers by monochromatic pieces — a survey of results and
  problems.
\newblock {\em Discrete Math.}, 339(7):1970--1977, 2016.

\bibitem{Gyarfas2006}
A.~Gyárfás, M.~Ruszinkó, G.~N. Sárközy, and E.~Szemerédi.
\newblock An improved bound for the monochromatic cycle partition number.
\newblock {\em J. Combin. Theory Ser. B}, 96(6):855--873, 2006.

\bibitem{Han2017}
J.~Han, A.~Lo, and N.~Sanhueza-Matamala.
\newblock Covering and tiling hypergraphs with tight cycles.
\newblock {\em Combin. Probab. Comput.}, 30(2):288--329, 2021.

\bibitem{Haxell2009}
P.~E. Haxell, T.~{\L}uczak, Y.~Peng, V.~R\"{o}dl, A.~Ruci\'{n}ski, and
  J.~Skokan.
\newblock The {R}amsey number for 3-uniform tight hypergraph cycles.
\newblock {\em Combin. Probab. Comput.}, 18(1-2):165--203, 2009.

\bibitem{Lang2020}
R.~Lang and N.~Sanhueza-Matamala.
\newblock Minimum degree conditions for tight {H}amilton cycles.
\newblock {\em J. Lond. Math. Soc.}, 105(4):2249--2323, 2022.

\bibitem{Luczak1999}
T.~{\L}uczak.
\newblock \({R(C_n,C_n,C_n)} \leq (4+o(1))n\).
\newblock {\em J. Combin. Theory Ser. B}, 75(2):174--187, 1999.

\bibitem{Luczak1998}
T.~{\L}uczak, V.~Rödl, and E.~Szemer{\'e}di.
\newblock Partitioning two-coloured complete graphs into two monochromatic
  cycles.
\newblock {\em Combin. Probab. Comput.}, 7(4):423–436, 1998.

\bibitem{Pokrovskiy2014}
A.~Pokrovskiy.
\newblock Partitioning edge-coloured complete graphs into monochromatic cycles
  and paths.
\newblock {\em J. Combin. Theory Ser. B}, 106:70--97, 2014.

\bibitem{Pokrovskiy2016}
A.~Pokrovskiy.
\newblock Partitioning a graph into a cycle and a sparse graph.
\newblock {\em Discrete Math.}, 346(1):113161, 2023.

\bibitem{Sarkozy2014}
G.~N. S\'ark{\"o}zy.
\newblock Improved monochromatic loose cycle partitions in hypergraphs.
\newblock {\em Discrete Math.}, 334:52--62, 2014.

\end{thebibliography}
\end{document}